\documentclass[11pt]{amsart}
\usepackage{amsfonts}

\newtheorem{theorem}{\bf Theorem}
\newtheorem{remark}{\bf Remark}[section]
\newtheorem{proposition}{Proposition}[section]

\newtheorem{corollary}{Corollary}[section]

\newtheorem{example}{Example}[section]
\newtheorem{definition}{Definition}

\usepackage{enumerate,mathrsfs}
\usepackage{amssymb,amsmath,graphicx,amsthm,mathtools,hyperref}

\usepackage{booktabs}

\theoremstyle{plain}
\usepackage{siunitx}

\usepackage{algorithm}
\usepackage{algorithmic,colortbl}

\makeindex
\usepackage[intoc]{nomencl} 

\makenomenclature

\usepackage{hyperref}
\hypersetup{nesting=true,debug=true,naturalnames=true}
\usepackage{graphicx,amssymb,upref}

\usepackage{subcaption}
\usepackage{enumerate,float}
\usepackage{fullpage}

\begin{document}

\title{Geometrical subordinated Poisson processes and its extensions}
\author[]{Neha Gupta}
\address{\emph{Department of Mathematics and Statistics,
		Indian Institute of Technology Kanpur, Kanpur 208016, India.}}
\email{nehagpt@iitk.ac.in}
\author[]{Aditya Maheshwari}
\address{\emph{Operations Management and Quantitative Techniques Area, Indian Institute of Management Indore, Indore 453556, India.}}
\email{adityam@iimidr.ac.in}
\author[]{Dheeraj Goyal}
\address{\emph{ Department of Mathematics and Statistics,
		Indian Institute of Technology Kanpur, Kanpur 208016, India.}}
\email{dheerajg@iitk.ac.in}

\begin{abstract}
In this paper, we study a generalized version of the Poisson-type process by time-changing it with the geometric counting process. Our work generalizes the work done by Meoli (2023) \cite{meoli2023some}. We defined  the geometric subordinated Poisson process (GSPP), the geometric subordinated compound Poisson process (GSCPP) and the geometric subordinated multiplicative Poisson process (GSMPP) by time-changing the subordinated Poisson process, subordinated compound Poisson process and subordinated multiplicative Poisson process with the geometric counting process, respectively. We derived several distributional properties and many special cases from the above-mentioned processes. We calculate the asymptotic behavior of the correlation structure. We have discussed applications of time-changed generalized compound Poisson in shock modelling. 
\end{abstract}
\maketitle
%\begin{enumerate}
 %   \item geometric subordinated Poisson process (GSFPP)
  %  \item geometric subordinated compound Poisson process (GSCPP)
  %  \item geometric subordinated multiplicative Poisson process (GSMPP)
%\end{enumerate}

 \section{Introduction}

%\cite{meoli2023some,DiCrescenzo2019,Cha2013}
%The geometric process modifies Poisson processes by introducing random time scales governed by a geometric random variable. It addresses limitations like independent increments and memorylessness, offering geometric distributed increments, overdispersion, and infinite inter-event expectations. This model is widely used in reliability, seismology, and environmental sciences, providing practical solutions for probabilities, moments, and crossing times.\\
A counting process is a stochastic model used to track the occurrence of discrete events over time, often applied in fields like finance, reliability, and epidemiology. The Poisson process is the most popular counting process which serves as the standard framework for modeling random arrivals in continuous time; while its is an extremely important count model but its reliance on exponentially distributed interarrival times limits its applicability in some practical situations.  More specifically, it assumes independence and stationarity, but real-world events often deviate due to clustering, time-dependent rates, or inter-event dependencies. Identifying these departures is crucial for accurate modeling in fields like epidemiology and finance, where event patterns are more complex. In \cite{Cha2013}, authors introduced a counting process called the geometric counting process, where increments of the counting process follow the geometric distribution. It is a special type of stochastic process that builds on the concept of mixed Poisson processes, where the event rate is random rather than fixed. It originated from studies on mixed Poisson models, particularly those reviewed by \cite{grandell2020mixed,rolski2009stochastic}, and became specific when the rate parameter was modelled using an exponential distribution. This led to a process with geometric distributed increments and dependent events, unlike the traditional Poisson process. \\\\
The geometric counting process was further studied by \cite{DiCrescenzo2019}, which explored applications in shock models, seismology and software reliability. %\textcolor{red}{reference}.
These recent developments led to increased interest among researchers, and they have started using the geometric process as a time-change. In \cite{meoli2023some}, the author explored some Poisson-based processes where they used a geometric counting process as a time-change. The time-changed processes,  sometimes called as a subordinated stochastic process and their first-exit times, have been widely studied in the literature. The most general among them is a non-decreasing L\'evy process \cite{sato}, called a subordinator. Several well-known examples of subordinators include the gamma process, lower incomplete gamma process, inverse Gaussian process, stable process, tempered stable process, geometric stable, and \textit{etc}. %Subordinators are fundamental in defining a new class of stochastic processes known as subordinated stochastic processes. 
In this direction, researchers have studied subordinated Poisson process (SPP), which is obtained by applying a time-change by an independent L\'evy subordinator or its right continuous inverse to a Poisson process. More specific and notable examples are fractional Poisson process (see \cite{lask}), the space-fractional Poisson process \cite{Orsingher2012}, the tempered space-fractional Poisson process \cite{gupta2020tempered} and the space-fractional negative binomial process \cite{lrd2016}, \textit{etc}. The use of generalized L\'evy subordinator and/or its inverse is well-known in literature (see  \cite{TCFPP-pub, toaldo2015convolution, kataria2022generalized, polito2016generalization, soni2024bivariate, gupta2023fractional, gupta2024tempered}). These processes have various applications in finance, statistical physics \cite{uchaikin2008fractional, laskin2009some}, anomalous diffusion modelling \cite{piryatinska2005models}, fractional partial differential equations, \textit{etc.}\\\\
The compound Poisson and multiplicative compound Poisson process were also studied after time-changing them with geometric subordinator. The above-mentioned processes extends the Poisson process by allowing random jump sizes, making it ideal for modeling count data with varying event impacts. We have studied their distributional properties. These processes are widely used in fields such as insurance \cite{dickson2016insurance}, evolutionary biology \cite{huelsenbeck2000compound}, and reliability \cite{wang2022explicit}, where extreme events of random size are significant. Particularly, in the domain of reliability theory, we apply these processes in study of shock models: extreme shock model and cumulative shock model; these shock models have huge applications and widely studied in the literature. Despite the vast applications, most of the literature consider Poisson processes (homogeneous/non-homogeneous Poisson processes and mixed Poisson processes) to model shock arrivals. However, Poisson process are unsuitable to model events which have heavy-tailed inter-arrival times and/or multiple occurrences. The process which we study in the paper is highly flexible to capture such phenomena as fractional Poisson process,
the space-fractional Poisson process, and the tempered space-fractional Poisson process are some special cases of it. Thus by considering such a flexible process to model shock arrivals may help in better prediction of the system lifetime under shock environment. As of our knowledge, we are first to study shock models with such a general process which take cares of multiple occurrences as well as heaviness between inter-arrival times.\\\\
In this paper, we study several generalized versions of Poisson-based processes time-changed by geometric counting process. We call them as the geometric subordinated Poisson process (GSPP) the geometric subordinated compound Poisson process (GSCPP) and the geometric subordinated multiplicative Poisson process (GSMPP) by time-changing the subordinated Poisson process, subordinated compound Poisson process and subordinated multiplicative Poisson process with the geometric counting process, respectively. We derive distributional properties, first-passage time and asymptotic behavior of the correlation for the GSPP and its special cases.  Similar results were derived for GSCPP and GSMPP. Finally, we provide the applications of the above mentioned processes in extreme shock models.\\\\
The paper is organized as follows. In Section \ref{sec:prelims}, we provide some preliminary definitions and results. Section \ref{sec:subordinated} introduces the time-changed subordinated Poisson process with geometric time, along with a detailed discussion of its distributional properties. In Sections \ref{sec:subcomp} and \ref{sec:multcomp}, we explore the generalized fractional compound Poisson process with a geometric random component and the generalized fractional multiplicative Poisson process with a geometric random component, respectively. Finally, in Section \ref{sec:appl}, we present applications in shock models for specific processes.

% In many practical scenarios, such as customer arrivals in service systems, system failures in reliability engineering, and transaction occurrences in financial markets, events are naturally observed at discrete time points. This necessitates the development of discrete-time alternatives to traditional counting processes.

%The modeling of count data is 
 % The geometric process proved useful in fields where the time between events can have infinite expectations, such as earthquake modeling in seismology and system reliability studies. Researchers like Cha and Finkelstein highlighted its benefits for modeling events that are positively dependent. Over time, it has been widely applied in areas such as software reliability, failure analysis, and shock models, making it a versatile tool for analyzing complex random phenomena.
\ifx
\textcolor{red}{
\begin{enumerate}
    %\item The Geometric counting process has dependent waiting times. \cite{meoli2023some,DiCrescenzo2019,Cha2013}
    %\begin{enumerate}
    %    \item It is a mixed Poisson process. \cite{Cha2013}
    %    \item Applications, relevance and novelty for the geometric counting process.
    %    \item \cite{rolski2009stochastic}
    %\end{enumerate}
    \item Subordination of the geometric counting process \begin{enumerate}
        \item Some background on Poisson process subordination (some background by iterated Poisson, (fractional) Poisson, Pearson diffusion)\cite{ascione2021time}, \cite{orsingher2012space} space fractional Poisson process, \cite{maheshwari2019fractional} \cite{Orsingher2012}(considered the composition of a homogeneous Poisson process to an independent Poisson process, both homogenous and inhomogeneous; to an independent fractional birth process, both linear and non linear, and to an independent Poisson field. ), \cite{di2015compound}(studied a compound Poisson process whose running time is an independent Poisson process.)\cite{gupta2023fractional}(fractional generalization of CPP)
        \item Relevance, application and novelty \cite{zarezadeh2018network}(studied the reliability and stochastic properties of an n-component network assuming that the failure of each component is due to a Geometric Counting Process.)
        \item Subordination geometric counting process, compound geometric and multiplicative geometric its applications and that it is a research gap which we are addressing/closing.
            \end{enumerate}
    \item  What we have done in this paper?
    \item Stressing on the applications of the proposed approach.
 %   \item Brief survey of results and structure of paper.
        \end{enumerate}}\fi

\section{Preliminaries}\label{sec:prelims}
In this section, we recall some relevant definitions and properties of L\'evy subordinator, space-fractional Poisson process (SFPP), tempered space-fractional Poisson process (TSFPP), which will be used in analyzing the some generalized counting processes with geometric random times.
\subsection{L\'evy Subordinator}
A L\'evy subordinator (hereafter referred to as the subordinator) $\{D^{f}(t)\}_{t\geq0}$ is a non-decreasing L\'evy process and its Laplace transform (LT) (see \cite[Section 1.3.2]{appm}) has the form
			\begin{equation*}%\label{subordinator-LT}
				\mathbb{E}[e^{-s D^{f}(t)}]=e^{-tf(s)},
				\;{\rm where}\; 
				f(s)=b s+\int_{0}^{\infty}(1-e^{-s x})\nu(dx),~b\geq0, s>0,
			\end{equation*}
			is the Bernst\'ein function (see \cite{Bernstein-book} for more details). 
			Here $b$ is the drift coefficient and $\nu$ is a non-negative L\'evy measure on positive half-line satisfying 
			\begin{equation*}
				\int_{0}^{\infty}(x\wedge 1)\nu(dx)<\infty~~{\rm and}~~\nu([0,\infty))= \infty,
			\end{equation*}
			which ensures that the sample paths of $\{D^{f}(t)\}_{t \geq 0}$ are almost surely $(a.s.)$  strictly increasing.
            \noindent We discuss some special cases of strictly increasing subordinators. The following subordinators with Laplace exponent denoted by $f(s)$ are very often used in literature.
\begin{equation}\label{Levy_exponent}
f(s) =
         \begin{cases}
                  s^{\alpha}, & $(stable subordinator)$;\\
                 \sum_{i=1}^{n}{c_{i}s^{\alpha_{i}}} &$(mixed stable  subordinator)$;\\
                  (s+\mu)^{\alpha}-\mu^{\alpha},&$(tempered stable  subordinator)$;\\
                 \sum_{i=1}^{n}{c_{i}((s+\mu_{i})^{\alpha_{i}}-{\mu_{i}}^{\alpha_{i}})}, &$(mixture of tempered stable  subordinator)$;\\
                  p\log(1+\frac{s}{\beta}), &$(gamma subordinator)$;\\
                  
                   \delta(\sqrt{2s+\gamma^2}-\gamma), & $(inverse Gaussian subordinator)$,
			\end{cases}
\end{equation}
where $0<\alpha<1$, $c_{i}\geq{0},\; \sum_{i=1}^{n}{c_{i}}=  1$, $\mu >0$, $p >0$, $\beta > 0$, $\gamma>0$, and  $\delta>0$.

\subsection{Geometric Counting Process (GCP)}
A geometric counting process (GCP) $\{G_\mu(t)\}_{t \geq 0}$ with parameter $\mu >0$
is class of the mixed Poisson process (see \cite{grandell2020mixed}), whose one-dimensional distribution is given by
\begin{align*}
P_k^{\mu}=\mathbb{P}[G_\mu(t)=k]=\int_{0}^{\infty}\mathbb{P}[N(t)=k]dU_{\mu}(\lambda), \; k \in 0,1,\ldots,
\end{align*}
where, $\{N(t)\}_{t \geq 0}$ is Poisson process with intensity $\lambda >0$ and $dU_{\mu}(\lambda)$ is the exponential distributed with mean $\mu >0$.\\
%\noindent \textbf{Property}
    The GCP $\{G_\mu(t)\}_{t \geq 0}$  satisfies the following properties for fixed $\mu >0$, i.e. (see \cite{rolski2009stochastic})
    \begin{enumerate}[(a)]\label{pmf_GCP}
       \item $G_{\mu}(0)=0$
        \item $\mathbb{P}[G_\mu(t+s)-G_\mu(t)=k]=\frac{1}{1+\mu s}\left(\frac{\mu s}{1+\mu s}\right)^k, \; \; s, t \geq 0.$
        \end{enumerate}
The above property says that it has stationary increment property and the GCP is geometric distributed with parameter $
\frac{\mu s}{1+\mu s}$ for every length of the time interval $s >0$.\\
The GCP possesses the dependent increments property. The dependence structure in the increments of this process is given by the positive upper orthant dependent increments property, i.e., for any arbitrary integer $m \geq 2$ and $0 < t_{1} < t_{2} < \dots < t_{m},$
$$\mathbb{P}(G_{\mu}(t + \Delta t_{i}) - G_{\mu}(t_{i}) > n_{i}, i = 1,2,\dots, m ) \geq \prod_{i=1}^{m} \mathbb{P}(G_{\mu}(t_{i} + \Delta t_{i}) - G_{\mu}(t_{i}) > n_{i}), $$
for all $n_{i}$, $i = 1,2,\dots, m$, where $t_{i} + \Delta t_{i} \leq t_{i+1},$ $i = 1,2,\dots, m-1$ (see Cha and Mercier~\cite{cm}). Such a dependency in the increments indicates that when the number of events in the past is higher then the likelihood of events in the future is also higher. This type of dependency is useful in many real-life applications; particularly, in reliability applications.   \\
 The LT of the $G_\mu(t)$ is given by
 \begin{align}\label{lt_gp}
\mathbb{E}[e^{-sG_\mu(t)}]=\frac{1}{1+\mu t(1-e^{-s})},\;\; \mathbb{R}(s)>0.
 \end{align}
 For a complex variable $y$, the geometric polynomials $w_n(y)$ of degree $n \in \mathbb{N}\cup \{0\}$ were introduced in Euler’s work \cite{boyadzhiev2007apostol,boyadzhiev2015power,}, are defined by
 $$
 w_n(y)=\sum_{k=0}^n { n \brace k} k!y^k,
 $$
 where 
 $$
 { n \brace k}=\frac{1}{k!}\sum_{i=0}^k(-1)^i{k \choose i}(k-i)^n,\;\; k=0,1,\ldots,n,
 $$
 are the Stirling numbers of the second kind.
 The geometric-like series
 \begin{align}\label{geo_series}
     \sum_{k=0}^\infty k^n y^k=\frac{1}{1-y}w_n\left(\frac{y}{1-y}\right),
 \end{align}
 when $|y| \leq 1$ and $n \in \mathbb{N}$.\\
 The mean, variance and the covariance of the geometric counting process $G_\mu(t)$ is given by
\begin{align}\label{mean_GP}
\mathbb{E}(G_\mu(t)) &= \mu t.\nonumber\\
\mathrm{Var}(G_\mu(t)) &= \mu t(1+\mu t).\\
{\rm Cov}[G_\mu(t),G_\mu(s)]&= \mu s(1+\mu t)\nonumber.
\end{align}
%\begin{corollary}\textcolor{red}{may be shifted to proof}
%		We can obtain the correlation function of $G_\mu(t)$ and $G_\mu(s)$. It exhibits SRD property, i.e.
%		$$
%		\lim\limits_{t\to\infty}\frac{{\rm Corr}[G_\mu(t), G_\mu(s)]}{t^{-1}}  \sim  \frac{s^{1/2}}{\lambda^{1/2}(1+\lambda s)^{1/2}}.
%		$$
		%Because $Y_{\alpha, \mu}(t)$ decays as $t^{-1/2}$.
%	\end{corollary} 
\subsection{Subordinated Poisson Process (SPP)}The process $\{N^f(t)\}_{ t > 0}$, can be viewed as subordinated Poisson processes where $\{D^f(t)\}_{t >0}$ are L\'evy subordinators associated with the Bernstein function $f$ and independent from the homogeneous Poisson process $\{N(t)\}_{t >0}$ with parameter $\lambda >0$ (see \cite{OrsToa-Berns}).\\
 The LT of the $N^f(t)=N(D^f(t))$ is given by
 $$\mathbb{E}[e^{-uN^f(t)}]= e^{-tf(\lambda(1-e^{-u}))}$$
 The mean, variance and covariance functions of the SPP are given by
\begin{align}\label{mean_spp}
\mathbb{E}[N^f(t)]&=\lambda \mathbb{E}[D^f(t)]\nonumber\\
 \mathrm{Var}(N^f(t))&=\lambda^2 \mathrm{Var}[D^f(t)]+\lambda \mathbb{E}[D^f(t)]
\end{align}
Let $\{N^f(t)\}_{t \geq 0}$ be a 
Orsingher and Toaldo \cite{OrsToa-Berns} studied the Poisson process subordinated with independent L\'evy subordinator.  
\noindent The pmf $P^{ f}_k(t)=\mathbb{P}[N^f(t)=k]$ of the $N^f(t)$ is given by
\begin{align}\label{pmf_shpp}
P^{f}_k(t)=\frac{(-1)^k}{k!}\frac{d^k}{du^k}e^{-tf(\lambda u)}|_{(u=1)}, \;\; k \geq 0.
\end{align}
The well known processes such as the space-fractional (SFPP) $\{N^{\alpha}(t)\}_{t\geq 0}$ and the tempered space-fractional Poisson process $\{N^{\alpha, \nu}(t)\}_{t\geq 0}$ (TSFPP) are obtained with time change in the Poisson process with independent $\alpha$-stable subordinator and tempered $\alpha$-stable subordinator, respectively.\\
The pmf of the $N^{\alpha}(t)$ is given by (see \cite{orsingher2012space})
\begin{align}\label{pmf_SFPP}
\mathbb{P}[N^{\alpha}(t)=k]=\frac{(-1)^k}{k!}\sum_{r=0}^{\infty}\frac{(-\lambda^\alpha t)^r}{r!}\frac{\Gamma(\alpha r+1)}{\Gamma(\alpha r+1-k)}
\end{align}
The pmf of the $N^{\alpha,\nu}(t)$ is given by (see \cite{gupta2020tempered})
\begin{align}\label{pmf_TSFPP}
\mathbb{P}[N^{\alpha, \nu}(t)=k]=\frac{(-1)^k}{k!}\frac{\lambda^k e^{\nu^\alpha t}}{(\lambda+\nu)^k}\sum_{m=0}^{\infty}\frac{(-t(\lambda+\nu))^m}{m!}\frac{\Gamma(\alpha m+1)}{\Gamma(\alpha m+1-k)}
\end{align}
The mean, variance and the covariance \cite{OrsToa-Berns} of the TSFPP is given by
\begin{align}\label{cov_TFSPP}
\mathbb{E}(N^{\alpha,\nu}(t)) &= \lambda \alpha \nu^{\alpha-1}t.\nonumber\\
\mathrm{Var}(N^{\alpha,\nu}(t)) &= \mathbb{E}(\lambda D_{\alpha,\mu}(t)) + \mathrm{Var}(\lambda D_{\alpha,\mu}(t)) = \lambda \alpha \nu^{\alpha-1}t + \lambda^2 \alpha(1-\alpha) \nu^{\alpha-2}t.\\
{\rm Cov}[N^{\alpha,\nu}(t),N^{\alpha,\nu}(s)] & = \lambda \alpha \nu^{\alpha-1}t + \lambda^2 \alpha(1-\alpha) \nu^{\alpha-2}s\nonumber.
\end{align}
\section{geometric subordinated Poisson processes}\label{sec:subordinated}
In this section, we study the generalized Poisson process time-changed by an independent geometric process and explore their distributional properties. We also study several special cases of the defined process.
\begin{definition}We define the geometric subordinated Poisson process (GSPP)  by time-changing the subordinated Poisson process (SPP)  $\{N^f(t)\}_{t \geq 0}$ with an independent geometric counting process (GCP)
$\{G_\mu(t)\}_{t \geq 0}$, such as 
%The time-changed subordinated Poisson process with geometric time $\{N^{f}_{\mu}(t)\}_{t \geq0}$ is defined as
\begin{align}\label{spp_g}
N^{f}_{\mu}(t):=N^f(G_\mu(t)), \;\; t \geq 0.
\end{align}
\end{definition}
The LT of the $N^{f}_{\mu}(t)$ is given by
\begin{align}\label{LT_GFPP_GP}
\mathbb{E}[e^{-uN^{f}_{\mu}(t)}]=\mathbb{E}[\mathbb{E}[e^{-uN^{f}_{\mu}(t)}|G_\mu(t)]] =\frac{1}{1+\mu t (1-e^{-f(\lambda (1-e^{-u}))})}.
\end{align}
\begin{proposition}\label{pr1}
      The probability mass function (pmf) $P^{ f, \mu}_k(t)=\mathbb{P}[N^{f}_{\mu}(t)=k]$ of the $N^{f}_{\mu}(t)$ is given by
    %  $$
     % P^{f,\mu}_k(t)=\frac{\lambda^k}{k!}\sum_{n=0}^{\infty}\frac{(-\lambda)^n}{n!}\mathbb{E}[D^f^{(n+k)}(G_\mu(t))]
     % $$
     % or \textcolor{red}{confusion}
$$
      P^{f, \mu}_k(t)=\frac{(-1)^k}{k!}\frac{d^k}{du^k}\left(\frac{1}{1+\mu t(1-e^{-f(\lambda u)})}\right)_{u=1}.
      $$
 \begin{proof}
       By using conditional arguments, we have
      \begin{align*}
          P^{f, \mu}_k(t) &=\mathbb{E}[\mathbb{P}[N^{f}_{\mu}(t)=k|G_{\mu}(t)]] \;\;\textit{putting the  pmf of SPP \eqref{pmf_shpp}},\\
          &= \sum_{n=0}^{\infty}\frac{1}{1+\mu t}\frac{(-1)^k}{k!}\frac{d^k}{du^k}e^{-nf(\lambda u)}|_{(u=1)}\left(\frac{\mu t}{1+\mu t}\right)^n,\\
          &= \frac{1}{1+\mu t}\frac{(-1)^k}{k!}\frac{d^k}{du^k}\left(\sum_{n=0}^{\infty}\left(\frac{e^{-f(\lambda u)} \mu t}{1+\mu t}\right)^n\right)_{u=1},
      \end{align*}
      By applying infinite sum of geometric series and cancel some term, we get the pmf.
 \end{proof}      
  \end{proposition}
 
  %\noindent In addition, we discuss some distributional properties of the process.
%		\begin{theorem}\label{m_v_c_GFCPP}
			
%			The mean, variance and covariance are given by:
%			\begin{enumerate}[(i)]
%				\item $ \mathbb{E}[N^{f}_{\mu}(t)] =\lambda \mathbb{E}[D^{f}(G_\mu(t)]$
%				\item ${\rm Var}[N^{f}_{\mu}(t)] = \lambda \mathbb{E}[D^{f}(G_\mu(t)]+\lambda^2{\rm Var}[D^{f}(G_\mu(t)]$
%				\item ${\rm Cov}[N^{f}_{\mu}(t),N^{f}_{\mu}(s)]  = \lambda \mathbb{E}[D^{f}(G_\mu(s)]+\lambda^2{\rm Cov}[D^{f}^2(G_\mu(t),D^{f}^2(G_\mu(s)].$
				% {\rm Cov}[Y_f(t),Y_f(s)] &= \lambda \mathbb {E}[E_f(s)]{\rm Var}[X_1]+( \mathbb{E}[X_1])^2{\rm Cov}[N_f(t),N_f(s)],\;\; s<t,
%			\end{enumerate}
%		\end{theorem}
  
% \begin{proof}
  %    $\mathbb{E}[N^{f}_{\mu}(t)]= \mathbb{E}[\mathbb{E}[N^{f}_{\mu}(t)|D^f(G_\mu(t))]]=\lambda \mathbb{E}[D^f(G_\mu(t))] $
  %    the varience 
 % \end{proof}
 \begin{proposition}\label{pr2}
    Let ${N^f_{\mu}(t)}_{t\geq 0}$ be the GSPP, as defined in \eqref{spp_g}, we have the following distributional equality in the compound form 
   $$ N^f_{\mu}(t) \stackrel{d}{:=} \sum_{i=1}^{G_\mu(t)}Y_i, \; \; t \geq 0,
$$ 
where $Y_i's$ are subordinated Poisson random variables are i.i.d. $N^{f}(1)$. 
  \end{proposition}
  \begin{proof}
     The LT of the \textit{pdf} of $\sum_{i=1}^{G_\mu(t)}Y_i$ is
			$$
			\mathbb{E}[e^{-u\sum_{i=1}^{G_\mu(t)}Y_i}] = \frac{1}{1+\mu t(1-M_{N^f(1)}(-u))}= \frac{1}{1+\mu t (1-e^{-f(\lambda (1-e^{-u})))})}.
			$$
			Comparing the above equation with the LT \eqref{LT_GFPP_GP} of the $\{N^{f}_{\mu}(t)\}_{t\geq 0}$, we get the desired result. 
  \end{proof}
Here, we consider the such kind of the SPP for which the moments $\mathbb{E}[(N^{f}(t))^r]<\infty$ for all $r >0$.  In addition, we discuss some distributional properties of the process.
		\begin{theorem}\label{m_v_c_GFCPP}
			
			The mean, variance and covariance is given by:
			\begin{enumerate}[(i)]
				\item $ \mathbb{E}[N^{f}_{\mu}(t)] =\lambda \mu t\mathbb{E}[D^{f}(1)]$
				\item ${\rm Var}[N^{f}_{\mu}(t)] =\lambda \mu t (\lambda\mathbb{E}[{D^f(1)}^2]+\mathbb{E}[D^f(1)])+ \mu^2 \lambda^2 t^2\mathbb{E}[D^f(1)]^2 $
				\item ${\rm Cov}[N^{f}_{\mu}(t),N^{f}_{\mu}(s)]  = \lambda \mu s (\lambda\mathbb{E}[{D^f(1)}^2]+\mathbb{E}[D^f(1)])+ \mu^2 \lambda^2 s t\mathbb{E}[D^f(1)]^2, \;\; 0<s <t.$
				% {\rm Cov}[Y_f(t),Y_f(s)] &= \lambda \mathbb {E}[E_f(s)]{\rm Var}[X_1]+( \mathbb{E}[X_1])^2{\rm Cov}[N_f(t),N_f(s)],\;\; s<t,
			\end{enumerate}
		\end{theorem}
  
 \begin{proof}
 Using \eqref{pr2}, we have that
      $\mathbb{E}[N^{f}_{\mu}(t)]= \mathbb{E}[N^f(1)]\mathbb{E}[G_\mu(t)]$.\\
      The variance of $\{N^f_{\mu}(t)\}_{t \geq 0}$ can be written as (see  \cite{LRD2014})
				$$
				{\rm Var}[N^f_{\mu}(t)] = {\rm Var}[N^f_{\mu}(t)]\mathbb{E}[G_{\mu}(t)]+\mathbb{E}[N^f(1)]^2{\rm Var}[G_{\mu}(t)].
				$$
                Next, we compute the ${\rm Cov}[N_{\mu}^f(t),N^f_{\mu}(s)]$ by using \eqref{pr2},  $s\leq t$,
				\begin{align*}
					{\rm Cov}[N_{\mu}^f(t),N^f_{\mu}(s)]&= \mathbb{E}\left[\sum_{k=1}^{\infty}Y_k^2 \mathbb{I}\{G_{\mu}(s)\geq k\}\right]+\mathbb{E}\left[ \sum_{i\neq j}\sum Y_i Y_k \mathbb{I}\{G_\mu(s)\geq k, G_\mu(t)\geq i\}\right]\\
					& \hspace{3cm}-(\mathbb{E}[Y_1])^2\mathbb{E}[G_\mu(s)]\mathbb{E}[G_\mu(t)]\\
					&= \mathbb{E}[Y_1^2]\sum_{k=1}^{\infty}\mathbb{P}(G_{\mu}(s)\geq k)\\
					&\hspace{3cm}+(\mathbb{E}[Y_1])^2\left[\sum_{i=1}^{\infty}\sum_{k=1}^{\infty}\mathbb{P}(G_{\mu}(s)\geq k, G_\mu(t)\geq i)-\mathbb{P}(G_\mu(s)\geq k)\right]\\
					&\hspace{3cm}-(\mathbb{E}[Y_1])^2\mathbb{E}[G_\mu(s)]\mathbb{E}[G_\mu(t)]\\
					&=  \mathbb{E}[Y_1^2] \mathbb{E}[G_\mu(s)]+(\mathbb{E}[Y_1])^2\mathbb{E}[G_\mu(s) G_\mu(t)]\\
					&\hspace{3cm}-(\mathbb{E}[Y_1])^2\mathbb{E}[G_\mu(s)]-(\mathbb{E}[Y_1])^2\mathbb{E}[G_\mu(s)]\mathbb{E}[G_\mu(t)]\\
					&= {\rm Var}[Y_1]\mathbb{E}[G_\mu(s)]+(\mathbb{E}[Y_1])^2{\rm Cov}[G_\mu(t),G_\mu(s)].
				\end{align*}
				where $Y_i's$ are subordinated Poisson distributed $N^f(1)$. Substituting equations \eqref{mean_spp} and \eqref{mean_GP} in the above equation, we get the desired result.\qedhere    
  \end{proof} 
  %\begin{example}(Fractional Negative Binomial process)
  %  Let $ D^f(t)$ be the gamma subordinator with Bernstein function $f(s)=p \log(1+\frac{s}{\beta})$ and LT is given by
%$$
%\mathbb{E}[e^{-u D^g (t)}]= \left(1+\frac{u}{\beta}\right)^{-pt},
%$$
%and 
%%Here, we define by time-changing the fractional negative binomial process by an independent geometric process, that is,
%$$N^{g, \mu}(t)= N^g(G_\mu(t))),$$
%where  $N^g(t)$  is fractional negative binomial process (see \cite{fnbpfp}).
%\begin{proposition}
 %   The pmf $P^{ g, \mu}_k(t)= \mathbb{P}[N^g\mu(t)=k]$ of the process $N^{g, \mu}(t)$ is given by
 %   $$
 %P^{ g,\mu}_k(t) = \frac{\lambda^k}{k!}\sum_{n=0}^{\infty} \frac{(-\lambda)^n}{n!}\sum_{j=0}^{\infty}\frac{\Gamma(\rho +pj)}{\beta^{\rho}\Gamma(pj)}\frac{1}{1+\mu t}\left(\frac{\mu t}{1+\mu t}\right)^j
 %$$
%\end{proposition}
 
  %\end{example}
  \noindent Next, we explore special cases of subordinated Poisson processes, including SFPP and TSFPP, and discuss their distributional properties.
\begin{example}(Geometric space fractional Poisson process) Let $\{N^{\alpha}_{\mu}(t)\}$
 be the time-changed in space-fractional Poisson process with GCP, is defined by time-changing the SFPP $\{N^{\alpha}(t)\}_{t \geq 0}$ by an independent GCP, that it
 $$
 N^{\alpha}_{\mu}(t)=N^\alpha(G_{\mu}(t)),\;,\, t \geq 0.
 $$
The LT of the $N^{\alpha}_{\mu}(t)$ is given by
\begin{align}\label{LT_SFPP_GP}
\mathbb{E}[e^{-uN^{\alpha}_{\mu}(t)}]=\mathbb{E}[\mathbb{E}[e^{-uN^{\alpha}_{\mu}(t)}|G_\mu(t)]] =\frac{1}{1+\mu t (1-e^{-\lambda^\alpha (1-e^{-u})^\alpha})}.
\end{align}     
    \end{example}  
  \begin{proposition}\label{pmf_sfpp_gcp}

  The pmf $\mathbb{P}[N^{\alpha}_{\mu}(t)=k]=P^{\alpha, \mu}_k(t)$ of the $N^{\alpha}_{\mu}(t)$ for any $k \in \mathbb{N}_{0}$ is given by 
  $$
  P^{\alpha, \mu}_k(t) = \sum_{r=0}^\infty {\alpha r \choose k} \frac{(-\lambda^\alpha)^r (-1)^k }{r!}w_r(\mu t),\;\; \lambda>0, \; \mu >0\;\; t \geq 0,
  $$
  where $w_{r}(\cdot)$ is the geometric polynomial of degree $r$.
  \begin{proof}
      By using conditional arguments and pmf  of SFPP, we have
        \begin{eqnarray}\label{deq1}
          P^{\alpha, \mu}_k(t) &=& \mathbb{E}[\mathbb{P}[N^{\alpha}_{\mu}(t)=k|G_{\mu}(t)]] \;\;\mbox{(substituting the  pmf of SFPP from \eqref{pmf_SFPP}})  \nonumber \\
          &=&   \frac{1}{1+\mu t}\sum_{r=0}^\infty {\alpha r \choose k} \frac{(-\lambda^\alpha)^r (-1)^k }{r!}\sum_{m=0}^\infty m^r \left(\frac{\mu t}{1+\mu t}\right)^m.
      \end{eqnarray}
      We use the geometric-like series property (see \cite{meoli2023some})
 \begin{align}\label{geo_series}
     \sum_{k=0}^\infty k^n y^k=\frac{1}{1-y}w_n\left(\frac{y}{1-y}\right),\; |y| \leq 1,\; n \in \mathbb{N}\cup\{0\}.
 \end{align}
      Using the above result for $y=\frac{\mu t}{1+\mu t}$, which is less than $1$, and rearranging the terms, completes the proof.
  \end{proof}
    \end{proposition}
 \noindent   We study the distribution of the first-passage time of the SFPP at geometric times at level k
  $$
  T_k=\inf\{s: N^{\alpha}(s)=k \},\;\; k \in \mathbb{N}.
  $$
  \begin{proposition}
      The distribution of the stopping time through an arbitrary level $k \in \mathbb{N}$ for
the process $N^{\alpha}_{\mu}(s)$ reads, for $s > 0$,
$$
\mathbb{P}[T_k \in ds]=ds\sum_{j=0}^{k-1}\sum_{r=0}^{\infty}{\alpha r \choose j}\frac{(-\lambda^\alpha)^r(-1)^{j+1}}{r!}\sum_{i=0}^{r}{r \choose i } \mu^i(i+1)!s^{i-1}.
$$
\begin{proof}That is,
    $$
    \frac{\mathbb{P}[T_k \in ds]}{ds}=\frac{d}{ds}\mathbb{P}[N^{\alpha}_{\mu}(s) \geq k]=
    \frac{d}{ds}\left[1-\sum_{j=0}^{k-1}\mathbb{P}[N^{\alpha}_{\mu}(s) = j]\right],
    $$
    substituting the eq \eqref{pmf_sfpp_gcp}, which completes the proof.
\end{proof}
  \end{proposition}   
  \begin{remark}
      Let $f(s)= s$, the process $\{N^{f}_{\mu}(t)\}_{t\geq0}$ reduces to  $N^{\mu}(t)=N(G_\mu(t))$ (see \cite{meoli2023some}). The asymptotic behaviour of the 
correlation of the process $N^{\mu}(t)$ and $N^{\mu}(s)$ for large time $t\geq 0$ is given by
  $$
  \lim\limits_{t\to\infty}\frac{{\rm Corr}[N^{\mu}(t), N^{\mu}(s)]}{t^{-1}}  \sim  \frac{s(1+\lambda)}{(\lambda \mu)^{1/2}\sqrt{s(1+\lambda+\lambda \mu s)}}.
  $$
  \end{remark}
 % \begin{proposition}
%The pmf $P^\alpha_k(t)$ satisfies the following governing equations
  %\begin{equation*}
%      \frac{d}{dt}P^\beta_k(t)=(-1)^k \sum_{r=0}^{\infty}
  %    \frac{(-\lambda ^\beta)^r}{r!}{\beta r \choose k}\frac{\mu}{1+\mu t}\left[\sum_{n=0}^{r-1} {r+1 \choose n} w_{n}(\mu t) +r w_r(\mu t)\right]
  %\end{equation*}
  %\end{proposition}
  \begin{example}(Geometric tempered space fractional Poisson process) The TSFPP $\{N^{\alpha,\nu}(t)\}_{t \geq 0}$ time-changed
by an independent GCP is defined as
$$
N^{\alpha,\nu}_{\mu}(t)=N^{\alpha,\nu}(G_\mu(t)),\;\; t \geq 0.
$$

\noindent Here, we discuss their distributional properties in details.
  \begin{proposition}
      The pmf $P^{\alpha, \nu,\mu}_k(t) = \mathbb{P}[N^{\alpha,\nu}_{\mu}(t)=k] $ is given by
$$
  P^{\alpha, \nu,\mu}_k(t)=\sum_{m=0}^{\infty}  \sum_{r=0}^\infty \nu^m \lambda^{\alpha r-m}{\alpha r \choose m}{\alpha r -m \choose k} \frac{ (-1)^k }{r!}\frac{1}{1+\mu t(1-e^{\nu ^\alpha})}w_r\left(\frac{e^{\nu^\alpha}\mu t}{1+\mu t(1-e^{\nu ^\alpha})}\right)
  $$
where $|e^{\nu^\alpha} \mu t| < 1+\mu t$. For $\nu=0$ reduces to .

 \begin{proof}
     This is follow the similar step of Prop \eqref{pmf_sfpp_gcp}.
 \end{proof}     
  \end{proposition}
  \begin{proposition}\label{mean_TSFPP_GP}
  We calculate the mean, variance and covariance of the $N^{\alpha, \nu}_{\mu}(t)$ is given by
			\begin{enumerate}[(i)]
				\item $ \mathbb{E}[N^{\alpha,\nu}_{\mu}(t)] =\lambda \alpha \nu^{\alpha-1}\mu t$
				\item ${\rm Var}[N^{\alpha,\nu}_{\mu}(t)] = \lambda \alpha \nu^{\alpha-1}\mu t + \lambda^2 \alpha(1-\alpha) \nu^{\alpha-2}\mu t+\lambda^2 \alpha^2\nu^{2(\alpha-1)}\mu t(1+\mu t)$
				\item ${\rm Cov}[N^{\alpha,\nu}_{\mu}(t),N^{\alpha,\nu}_{\mu}(s)]  = \alpha \lambda \nu^{\alpha-2}(\lambda(1-\alpha)+\nu)\mu s+\alpha^2\lambda^2\nu^{2(\alpha-1)}(\mu s+2 \mu^2 st).$
				% {\rm Cov}[Y_f(t),Y_f(s)] &= \lambda \mathbb {E}[E_f(s)]{\rm Var}[X_1]+( \mathbb{E}[X_1])^2{\rm Cov}[N_f(t),N_f(s)],\;\; s<t,
			\end{enumerate}
		\end{proposition}
 \begin{proof} 
 Using Theorem \eqref{m_v_c_GFCPP} and equation \eqref{cov_TFSPP}, we get the results.
 %Using eq \eqref{cov_TFSPP}, we get
 %\begin{align}\label{mean_TSFPP_Gp}
% \mathbb{E}[N^{\alpha,\nu}_{\mu}(t)]=\mathbb{E}[\mathbb{E}[N^{\alpha, \nu\mu}(t)|G_\mu (t)]]=\lambda \alpha \nu^{\alpha-1}\mathbb{E}[G_\mu (t)]=\lambda \alpha \nu^{\alpha-1}\mu t
% \end{align}
 %Further, form the eq \eqref{vari_TFSPP}
 %\begin{align*}
% \mathrm{Var}(N^{\alpha,\nu,\mu}(t)) = &\mathbb{E}[\mathrm{Var}(N^{\alpha,\nu}_{\mu}[G_{\mu}(t))|G_{\mu}(t)]] +  \mathrm{Var}[\mathbb{E}[N^{\alpha,\nu,\mu}(G_{\mu}(t))|G_{\mu}(t)]]\\
 %=& \lambda \alpha \nu^{\alpha-1}\mathbb{E}[G_\mu (t)] + \lambda^2 \alpha(1-\alpha) \nu^{\alpha-2}\mathbb{E}[G_\mu (t)]+\lambda^2 \alpha^2\nu^{2(\alpha-1)}\mathrm{Var}[G_\mu (t)]\\
% &= \lambda \alpha \nu^{\alpha-1}\mu t + \lambda^2 \alpha(1-\alpha) \nu^{\alpha-2}\mu t+\lambda^2 \alpha^2\nu^{2(\alpha-1)}\mu t(1+\mu t).
% \end{align*}
 %From eq \eqref{vari_TFSPP} and \eqref{cov_TFSPP}
% \begin{align}
 %\mathbb{E}[N^{\alpha,\nu}_{\mu}(t) N^{\alpha,\nu}_{\mu}(s)] &= \mathbb{E}[\mathbb{E}[N^{\alpha, \nu}(G_\mu (t)) N^{\alpha, \nu}(G_\mu(s))|(G_\mu(t), G_\mu(s))]]\nonumber\\
 %& = \alpha \lambda \nu^{\alpha-2}(\lambda(1-\alpha)+\nu)\mu s+\alpha^2 \lambda^2 \nu^{2(\alpha-1)}\mathbb{E}[G_\mu(t) G_\mu(s)]\nonumber,
 %\end{align}
 %Using the above expression and \eqref{mean_TSFPP_Gp}, Part (iii) follows.
 \end{proof}
 \begin{corollary}
     The asymptotic behaviour of the correlation of the process $N^{\alpha,\nu}_{\mu}(t)$ and $N^{\alpha,\nu}_{\mu}(s)$ is given by
  $$
  \lim\limits_{t\to\infty}\frac{{\rm Corr}[N^{\alpha,\nu}_{\mu}(t), N^{\alpha,\nu}_{\mu}(s)]}{t^{-1}}  \sim  \frac{s(1+\nu^{-1} \lambda(1-\alpha)+\alpha \lambda \nu^{\alpha-1}s)}{(\lambda \alpha \mu)^{1/2}\sqrt{ \nu^{\alpha-1} s + \lambda (1-\alpha) \nu^{\alpha-2} t+\lambda\alpha\nu^{2(\alpha-1)}s(1+\mu s)}}.
  $$
 \end{corollary}
 \noindent The index of dispersion for a counting process $Z(t)$ is defined by (see \cite{cox-lewis})
 $$
 I(t)=\frac{\mathrm{Var}(X
(t))}{\mathbb{E}[X(t)]}
 $$
 The stochastic process $X(t)$
 is said to be overdispersed if $I(t) > 1$
 for all (see \cite{BegClau14}). From Prop. \eqref{mean_TSFPP_GP}, we have that
  The index of dispersion is greater than $1$ for all $t \geq 0$, such that
  $$
  I(t)=\frac{\mathrm{Var}(N^{\alpha,\nu,\mu}
(t))}{\mathbb{E}[N^{\alpha,\nu}_{\mu}(t)]}=1+\lambda(1-\alpha)\nu^{-1}+\lambda \alpha \nu^{\alpha-1}(1+\mu t)>1,
  $$
this process exhibits overdispersion.
 \end{example}
 \section{Geometric subordinated compound and multiplicative Poisson processes}
 In this section, we will study two extensions of the Poisson processes, namely, the compound Poisson process and the multiplicative compound Poisson process. 
\subsection{Geometric subordinated compound Poisson process (GSCPP)}\label{sec:subcomp}
\begin{definition}
			 Let
			$X_i,\;i=1,2,\ldots,$ be the iid jumps with common
			distribution $F_X$ and let $N^f(t)$ be subordinated Poisson process with L\'evy subordinator. The process defined by
			\begin{equation}
				Y^{f}(t) := \sum_{i=1}^{N_{f}(t)} X_i,t\geq0,
			\end{equation}
			is called the generalized compound Poisson process (GCPP).
		\end{definition}
\noindent The LT of the $Y^{f}(t)$ is given by
\begin{align}\label{mgf_cpp}
    \mathbb{E}[e^{-uY^f(t)}]=e^{-t f(\lambda(1-M_X(-u)))},
\end{align}
where $M_X(u)$ is the LT of the $X$.

\begin{proposition}\label{mean_cpp}
The mean and variance of the process $Y^{f}(t)$ is given by:
\begin{enumerate}[(i)]\label{mean_cpp}
				\item $ \mathbb{E}[Y^{f}(t)] =\lambda \mathbb{E}[D^{f}(t)]\mathbb{E}[X_1]$
				\item ${\rm Var}[Y^{f}(t)] = \lambda \mathbb{E}[D^f(t)]\mathbb{E}[{X}^2]+\mathbb{E}[X]^2\lambda^2 {\rm Var}[D^f(t)] $
				\item ${\rm Cov}[Y^f(t),Y^f(s)]  = \lambda \mathbb{E}[X_1^2]\mathbb{E}[D^f(s)]+\lambda^2(\mathbb{E}[X_1])^2{\rm Cov}[D^f(t),D^f(s)]$.
			\end{enumerate}
 \begin{proof} These results follow on similar lines as \cite{ gupta2023fractional}, and therefore omitted.
 \end{proof}   
\end{proposition}
\noindent The cumulative distribution function (cdf) $W_{Y}(y, t)=\mathbb{P}[Y^f(t) \leq y] $ of $Y^f(t)$, for $y \in \mathbb{R}$ and $t >0$, is given by
$$
W_{Y}(y, t)= \sum_{m=0}^{\infty} P^f_m(t)F_Y^{(m)}(y).
$$
For $m \in \mathbb{N}$, the m-fold convolution of $F_Y(.)$ and $F^{(0)}_Y(y)=\mathbb{I}_{y \geq 0}.$

\begin{definition}
Consider the GCPP $\{Y^f(t)\}_{t \geq 0}$ time-changed with an independent GCP $\{G_\mu(t)\}_{t \geq0}$, defined as follows
$$
Y^f_{\mu}(t)=Y^f(G_{\mu}(t)), \; \; t \geq 0.
$$
The process $\{Y^f_{\mu}(t)\}_{t \geq 0}$ is referred to as the geometric subordinated compound Poisson process (GSCPP)
\end{definition}
\begin{proposition}
The LT of the $Y^f_{\mu}(t)$ is given by
$$
\mathbb{E}[e^{-uY^f_\mu(t)}]= \frac{1}{1+\mu t(1-e^{-f(\lambda(1-M_X(-u))})}, u>0.
$$
\begin{proof}From the eq. \eqref{mgf_cpp} and \eqref{lt_gp}, it follows 
 $$  
        \mathbb{E}[e^{-u Y^f_{\mu}(t)}]=\mathbb{E}[\mathbb{E}[e^{-uY^f(G_\mu(t))}|G_\mu(t)]]=\mathbb{E}[e^{-G_\mu(t)f(\lambda(1-M_X(-u))}].\qedhere
        $$
\end{proof}
\end{proposition}
\begin{proposition}
    We have the following equality in distribution for the GSCPP
$$
Y^f_{\mu}(t) \stackrel{d}{=} Z^\mu(t)=\sum_{i=0}^{G_\mu(t)}\chi_i,
$$
where $\chi_0=0$ and $\chi_1, \chi_2,\ldots$ are iid random variables distributed as $Y^f(1)$.
\begin{proof}
    It follows a similar approach as in \eqref{pr2}.
\end{proof}
\end{proposition}

\noindent Further, we discuss the mean and variance of the process GSCPP $\{Y_\mu^f(t)\}_{t\geq 0}$.
\begin{theorem}\label{m_v_c_gCPP_gp}
        For $0 <s \leq t <\infty$, The
distributional properties of the $\{Y^{f}_\mu(t)\}_{t \geq0}$ are as follows:
			
			\begin{enumerate}[(i)]
				\item $ \mathbb{E}[Y^{f}_{\mu}(t)] =\lambda \mu t \mathbb{E}[D^{f}(1)]\mathbb{E}[X_1]$
				\item ${\rm Var}[Y^{f}_{\mu}(t)] = \lambda \mu t [\mathbb{E}[D^{f}(1)]{\rm Var}[X_1]+\mathbb{E}[X^2](\lambda {\rm Var}[D^f(1)] +\mathbb{E}[D^f(1)])]+\lambda^2 \mu t \mathbb{E}[D^f(1)]^2\mathbb{E}[X_1]^2(1+\mu t)$
                \item ${\rm Cov}[Y^{f}_{\mu}(t),Y^{f}_{\mu}(s)]  = \mu s(\lambda \mathbb{E}[X_1^2]\mathbb{E}[D^f(1)]+\mathbb{E}[X_1^2]\lambda^2 {\rm Var}[D^f(1)] )+\lambda^2\mathbb{E}[D^f(1)]^2\mathbb{E}[X_1]^2 \mu s(1+\mu t)$
				% {\rm Cov}[Y_f(t),Y_f(s)] &= \lambda \mathbb {E}[E_f(s)]{\rm Var}[X_1]+( \mathbb{E}[X_1])^2{\rm Cov}[N_f(t),N_f(s)],\;\; s<t,
			\end{enumerate}
		\end{theorem}
        \begin{proof}
            Using the conditional argument and independence of $G_\mu(t)$ and $\chi_i$, we have that 
				$$
				\mathbb{E}[Y^{f}_{\mu}(t)]= \mathbb{E}[G_\mu(t)]\mathbb{E}[\chi_1]=\lambda \mu t\mathbb{E}[X_1]\mathbb{E}[D^f(t)].
				$$
				From eq \eqref{mean_cpp} and \eqref{mean_GP},the variance of $\{Y_f(t)\}_{t \geq 0}$ can be written as (see   \cite{LRD2014})
				$$
				{\rm Var}[Y^{f}_{\mu}(t)] = {\rm Var}[\chi_1]\mathbb{E}[G_\mu (t)]+\mathbb{E}[\chi_1]^2{\rm Var}[G_\mu(t)].
				$$
        Next, we can directly compute the ${\rm Cov}[Y^{f}_{\mu}(t),Y^{f}_{\mu}(s)]$,  $s\leq t$ by using the eq \eqref{mean_cpp}, such as
        $$
        {\rm Cov}[Y^{f}_{\mu}(t),Y^{f}_{\mu}(s)]={\rm Var}[Y^f(1)]\mathbb{E}[G_\mu(s)]+\mathbb{E}[Y^f(1)]^2{\rm Cov}[G_\mu(t),G_\mu(s)].
        $$
Now, substituting equations \eqref{mean_cpp} and \eqref{mean_GP} in the above expression, which completes the proof.
        \end{proof}
%\begin{corollary}
%    The process $Y^{f}_{\mu}(t)$ exhibit
%\end{corollary}
\begin{proposition}
    The cumulative distribution function (cdf) $W_Y^{\mu, f}(y, t):= \mathbb{P}[Y^f_{\mu}(t) \leq y]$ of the $Y^f_{\mu}(t)$ is given by
   % $$
  %  W_Y^\mu(y, t)=\frac{1}{1+\mu t(1-e^{-f(\lambda)})}\mathbb{I}_{y \geq 0}+\sum_{m=1}^{\infty} \frac{(-1)^m}{m!}F_Y^{(m)}(y)(y)\sum_{l=0}^{\infty}\frac{(-1)^lw_{l}(\mu t)}{l!} \frac{d^m}{du^m}(f(\lambda u)^l)_{u=1}
   % $$
    $$
    W_Y^{\mu}(y, t)=\frac{1}{1+\mu t(1-e^{-f(\lambda^{\alpha}}))}\mathbb{I}_{y \geq 0}+\sum_{m=1}^{\infty} \frac{(-1)^m}{m!}F_Y^{(m)}(y)\frac{d^m}{du^m}\left(\frac{1}{1+\mu t(1-e^{-f(\lambda u})}\right)_{u=1}
    $$
    \begin{proof}By using conditional arguments:
    \begin{align*}
    W_Y^{\mu,f}(y, t)&= \sum_{n=0}^{\infty}\mathbb{P}[G_{\mu}(t)=n]W_Y^f(y, n)\\
    &= \sum_{n=0}^{\infty}[G_{\mu}(t)=n]\left[\mathbb{I}_{\{y \geq 0\}}e^{-nf(\lambda)}+\sum_{m=1}^{\infty}P_m^f(n) F^{m}_Y(y)\right]\\
    &= \frac{1}{1+\mu t (1-e^{-f(\lambda)})}\mathbb{I}_{\{y \geq 0\}}+\sum_{n=0}^{\infty}\left(\frac{\mu t}{1+\mu t}\right)^n\frac{1}{1+\mu t}\sum_{m=1}^{\infty}F^{m}_Y(y)\frac{(-1)^m}{m!}\frac{d^m}{du^m}\left(e^{-nf(\lambda u)}\right)_{u=1}\\
    &= \frac{1}{1+\mu t (1-e^{-f(\lambda)})}\mathbb{I}_{\{y \geq 0\}}+\frac{1}{1+\mu t}\sum_{m=1}^{\infty}\frac{(-1)^m}{m!}F^{m}_Y(y)\frac{d^m}{du^m}\left(\sum_{n=0}^{\infty}\left(\frac{\mu t e^{-f(\lambda u)}}{1+\mu t}\right)^n \right)_{u=1}\\
     &= \frac{1}{1+\mu t (1-e^{-f(\lambda)})}\mathbb{I}_{\{y \geq 0\}}+\sum_{m=1}^{\infty}\frac{(-1)^m}{m!}F^{m}_Y(y)\frac{d^m}{du^m}\left(\frac{1}{1+\mu t(1-e^{-f(\lambda u)})}\right)_{u=1}.\qedhere
     \end{align*}
    \end{proof}
    \end{proposition}
   
%When $N^f(t)$ is SFPP then the pmf 
%$$
%P_k^{\alpha}(t)= \frac{(-1)^k}{k!}\frac{d^k}{du^k}(e^{-t \lambda^\alpha u^\alpha})_{u=1}.
%$$
%The cdf $W_Y^{\mu, \alpha}(y, t)= \mathbb{P}[Y^f_{\mu,\alpha}(t) \leq y]$ of the $Y^f_{\mu, \alpha}(t)$ is given by
 %   $$
%    W_Y^{\mu, \alpha}(y, t)=\frac{1}{1+\mu t(1-e^{-\lambda^{\alpha}})}\mathbb{I}_{\{y \geq 0\}}+\sum_{m=1}^{\infty} \frac{(-1)^m}{m!}F_Y^{(m)}(y)\sum_{r=0}^{\infty}{\alpha r \choose m }\frac{(-\lambda^\alpha)^r}{r!} w_{r}(\mu t)
 %   $$
%$$
 %   W_Y^{\mu, \alpha}(y, t)=\frac{1}{1+\mu t(1-e^{-\lambda^{\alpha}})}\mathbb{I}_{y \geq 0}+\sum_{m=1}^{\infty} \frac{(-1)^m}{m!}F_Y^{(m)}(y)(y)\frac{d^m}{du^m}\left(\frac{1}{1+\mu t(1-e^{-\lambda^\alpha u^\alpha})}\right)_{u=1}
 %   $$
    \begin{corollary}
        If $Y_i's$ are absolutely continuous with pdf $f_{Y}(\cdot),$ then 
        \begin{itemize}
            \item the continuous part of $f_y(\cdot)$ is given by $$g_{Y}^{f, \mu}(y, t)= \sum_{m=1}^{\infty} \frac{(-1)^m}{m!}F_Y^{(m)}(y)(y)\frac{d^m}{du^m}\left(\frac{1}{1+\mu t(1-e^{-f(\lambda u})}\right)_{u=1}
         $$
         \item a discrete part $f_y(\cdot)$ is given by 
         $$
         \mathbb{P}[Y^{f}(t)=0]=\frac{1}{1+\mu t(1-e^{-f(\lambda^{\alpha}}))}
         $$
        \end{itemize}
    \end{corollary}
    \begin{corollary}
        If $Y_n's$ are discrete and integer-valued random variables
        $$
        h_k^{*m}:= \mathbb{P}[Y_1+Y_2+\cdots+Y_m=k], k\in z,
        $$
        then 
        $$
        \mathbb{P}[Y^f_\mu(t)=k]=\frac{1}{1+\mu t(1-e^{-f(\lambda^{\alpha})})}\mathbb{I}_{k=0}+\sum_{m=1}^{\infty} \frac{(-1)^m}{m!}h_k^{*m}\frac{d^m}{du^m}\left(\frac{1}{1+\mu t(1-e^{-f(\lambda u)}}\right)_{u=1}.
        $$
%\subsection{Bernoulli jumps}      
%$\mathbb{P}[Y=k]=p^{k}(1-p)^{1-k}$ for $k \in (0,1)$.
%Now the $m$-fold convolution $Y$ is
%$$h^{*m}_k={m \choose k}p^k (1-p)^{m-k}$$
        
   \end{corollary}
    \begin{corollary}
        If $Y_n's$ are discrete and integer valued random variables
        $$
        h_k^{*m}:= \mathbb{P}[Y_1+Y_2+\cdots+Y_m=k], k\in z
        $$
        then 
        $$
        \mathbb{P}[Y^f_\mu(t)=k]=\frac{1}{1+\mu t(1-e^{-f(\lambda^{\alpha})})}\mathbb{I}_{k=0}+\sum_{m=1}^{\infty} \frac{(-1)^m}{m!}h_k^{*m}\frac{d^m}{du^m}\left(\frac{1}{1+\mu t(1-e^{-f(\lambda u)}}\right)_{u=1}
        $$
%\subsection{Bernoulli jumps}      
%$\mathbb{P}[Y=k]=p^{k}(1-p)^{1-k}$ for $k \in (0,1)$.
%Now the $m$-fold convolution $Y$ is
%$$h^{*m}_k={m \choose k}p^k (1-p)^{m-k}$$
        
    \end{corollary}
   %  \begin{proposition}
    %    The CDF $W_{Y}^{\alpha,\mu}(y,t)=\mathbb{P}[Y^\alpha_\mu(t) \leq y]$ of the space fractional compound Poisson process with an independent geometric time $Y^\alpha_\mu(t)$ is
    %    \begin{align}
     %       W_{Y}^{\alpha,\mu}(y,t)=\frac{1}{1+\mu t (1-e^{-\lambda^\alpha})}\mathbb{I}_{\{y \geq 0\}}+ \sum_{m=1}^{\infty}\frac{(-1)^m}{m!}F^{m}_Y(y)\sum_{r=0}^\infty {\alpha r \choose m} \frac{(-\lambda^\alpha)^r}{r!}w_r(\mu t)
    %    \end{align}
   % \end{proposition}
   % \begin{proof}By using conditional arguments:
   % \begin{align*}
   % W_Y^{\mu,\alpha}(y, t)&= \sum_{n=0}^{\infty}\mathbb{P}[G_{\mu}(t)=n]W_Y^\alpha(y, n)\\
   % &= \sum_{n=0}^{\infty}[G_{\mu}(t)=n]\left[\mathbb{I}_{\{y \geq 0\}}e^{-n\lambda^\alpha}+\sum_{m=1}^{\infty}P_m^{\alpha}(n) F^{(m)}_Y(y)\right]\\
   % &= \frac{1}{1+\mu t (1-e^{-\lambda^{\alpha}})}\mathbb{I}_{\{y \geq 0\}}+\frac{1}{1+\mu t}\sum_{m=1}^{\infty}F^{m}_Y(y)\frac{(-1)^m}{m!}\sum_{r=0}^{\infty}{\alpha r \choose m}\frac{(-\lambda^\alpha)^r}{r!}\sum_{n=0}^{\infty}\left(\frac{\mu t}{1+\mu t}\right)^n n^r,
    % \end{align*}
    % From \eqref{geo_series} with $y=\frac{\mu t}{1+\mu t}$, we get the cdf.
    % \end{proof}
    \subsection{Geometric subordinated multiplicative Poisson process (GSMPP)}\label{sec:multcomp}
    In this section, we introduce the generalized multiplicative CPP and study their properties.
It is to note that the  multiplicative CPP $\{N_{\pi}(t)\}_{t \geq 0}$ was introduced by Orsingher and Polito (2012) \cite{Orsingher2012}.

\begin{definition}
    Let $X_i's$ be iid random variables which are independent of subordinated Poisson process $N^f(t)$. Then the process $N_{\pi}^{f}(t)$ is defined as
        \begin{align}
			N_{\pi}^{f}(t)= \Pi_{i=1}^{N_{f}t)}X_i,, t\geq 0,
		\end{align}
		is called generalized multiplicative  CPP.
	\end{definition}
    \noindent The cdf $Q_{\pi}^{f}$ of the process $N_{\pi}^{f}(t)$ is given by
\begin{align}\label{cdf_multi_gen}
    Q_{\pi}^{f}(y;t) &= \sum_{m=0}^{\infty} \mathbb{P}[N^{f}(t)=m]\mathbb{P}[V_m \leq y]\nonumber\\
    &= e^{-f(\lambda)t}\mathbb{I}_{\{y \geq 1\}}+ \sum_{m=1}^{\infty} \mathbb{P}[N^{f}(t)=m]\mathbb{P}[V_m \leq y],\;\; y\in \mathbb{R}, \;\;t\geq0,
\end{align}
where $V_m=X_1X_2\ldots X_m$.
\begin{definition}
    We define the generalized multiplicative CPP by time-changing the $N_\pi^{f}(t)$ by an independent  GCP $G_\mu (t)$. It is given by
    $$
    Y_{\pi}^{f, \mu}(t):= N_\pi^{f}(G_\mu(t)),\; \mu >0.
    $$
We call it the generalized subordinated multiplicative Poisson process (GSMPP).
    \end{definition}
    \begin{proposition}\label{pr4.5}
        The cdf $Q_{\pi}^{f, \mu}(t)$ of the $Y_{\pi}^{f, \mu}(t)$ is given by
        \begin{align*}
        Q_{\pi}^{f,\mu}(y, t)=\frac{1}{1+\mu t(1-e^{-f(\lambda)})}\mathbb{I}_{y \geq 0}+\sum_{m=1}^{\infty} \frac{(-1)^m}{m!}F_Y^{(m)}(y)\frac{d^m}{du^m}\left(\frac{1}{1+\mu t(1-e^{-f(\lambda u})}\right)_{u=1}
    \end{align*}
    \begin{proof}
         From the equation \eqref{cdf_multi_gen}, we have that
        \begin{align*}
           Q_{\pi}^{f, \mu}(y; t)= &\sum_{n=0}^{\infty} \mathbb{P}[G_\mu(t)=n]Q_{\pi}^{f}(y; n)\\
           =&  \sum_{n=0}^{\infty} \mathbb{P}[G_\mu(t)=n]\Bigg[e^{-f(\lambda)t}\mathbb{I}_{\{y \geq 1\}}+ \sum_{m=1}^{\infty} \mathbb{P}[N^{f}(t)=m]\mathbb{P}[V_m \leq y]\Bigg],
        \end{align*}
substitute the pmf of the $G_\mu (t)$ and $N^{f}(t)$ form the equations \eqref{pmf_GCP} and \eqref{pmf_shpp} respectively, which completes the proof.
    \end{proof}
    \end{proposition}
\noindent Next, we calculate the Mellin transform of the $N_{\pi}^f(t)$ 
$$
\mathbb{E}[N_{\pi}^f(t)]^{\beta-1}=\sum_{k=0}^{\infty}\frac{(-1)^k}{k!}\mathbb{E}[X^{\beta-1}]^k \frac{d^k}{du^k}\left(e^{-t f(\lambda u)}\right)_{u=1}
$$
By using the previous result, we evaluated the Mellin transform of $Q_{\pi}^{f, \mu}(t)$ as folllows
\begin{align*}
\mathbb{E}[Y_{\pi}^{f, \mu}(t)]^{\beta-1}=& \frac{1}{1+\mu t}\sum_{n=0}^{\infty}\left(\frac{\mu t}{1+\mu t}\right)^{n}\mathbb{E}[N^{f}_\pi(n)]^{\beta-1}\\
&=\frac{1}{1+\mu t}\sum_{k=0}^{\infty}\frac{(-1)^k}{k!}\mathbb{E}[X^{\beta-1}]^k \frac{d^k}{du^k}\left(\sum_{n=0}^{\infty}\left(\frac{\mu t e^{-f(\lambda u)}}{1+\mu t}\right)^{n}\right)_{u=1}\\
&= \sum_{k=0}^{\infty}\frac{(-1)^k}{k!}\mathbb{E}[X^{\beta-1}]^k \frac{d^k}{du^k}\left(\frac{1}{1+\mu t(1-e^{-f(\lambda u)})}\right)_{u=1}
\end{align*}
Put $\beta=2$, then mean of the process is given 

\begin{eqnarray}\label{deqmcp}
\mathbb{E}[Y_{\pi}^{f, \mu}(t)] &=& \sum_{k=0}^{\infty}\frac{(-1)^k}{k!}\mathbb{E}[X]^k \frac{d^k}{du^k}\left(\frac{1}{1+\mu t(1-e^{-f(\lambda u)})}\right)_{u=1}.
\end{eqnarray}
We next study a special case of the GSMPP by taking the example of the tempered stable subordinator as a special case.
\begin{definition}
		Consider the subordinator associated with tempered stable Bernst\'ein function \eqref{Levy_exponent} $f(s)= (\nu+s)^\alpha-\nu^\alpha,\; \alpha \in (0,1], \nu>0$, denoted by $\{D_{\alpha,\nu}(t)\}_{t\geq0}$. Let $N_{\alpha,\nu}(t):=N(D_{\alpha,\nu}(t)),t\geq 0$ be the tempered space-fractional Poisson process (TSFPP) \cite{gupta2020tempered}. The process $N_{\pi}^{\alpha, \nu}(t)$ is defined by 
		\begin{align}
			N_{\pi}^{\alpha, \nu}(t)= \Pi_{i=1}^{N^{\alpha, \nu}(t)}X_i,, t\geq 0,
		\end{align}
		is called tempered multiplicative fractional CPP with  $X_i,\;i=1,2,\ldots,$ be the iid random variables independent of $N^{\alpha, \nu}(t)$.
	\end{definition}
\noindent The cdf $Q_{\pi}^{\alpha, \nu}$ of the process $N_{\pi}^{\alpha, \nu}(t)$ is given by
\begin{align}\label{cdf_multi_tem}
    Q_{\pi}^{\alpha, \nu}(y;t) &= \sum_{m=0}^{\infty} \mathbb{P}[N^{\alpha, \nu}(t)=m]\mathbb{P}[V_m \leq y]\nonumber\\
    &= e^{((\lambda+\nu)^{\alpha}-\nu^\alpha)t}\mathbb{I}_{\{y \geq 1\}}+ \sum_{m=1}^{\infty} \mathbb{P}[N^{\alpha, \nu}(t)=m]\mathbb{P}[V_m \leq y],
\end{align}
where $V_m=X_1X_2\ldots X_m$.
\begin{definition}
    We define new time-changed stochastic process,  characterized  by time evolving in a multiplicative CPP $N_\pi^{\alpha, \nu}(t)$ with GCP $G_\mu (t)$, such as
    $$
    Y_{\pi}^{\alpha, \nu, \mu}(t)= N_\pi^{\alpha, \nu}(G_\mu(t)),\; \nu, \mu >0, \alpha \in (0,1).
    $$
    \end{definition}
    \begin{proposition}\label{cdf_tem_multi}
        The cdf $Q_{\pi}^{\alpha, \nu, \mu}(t)= \mathbb{P}[Y_{\pi}^{\alpha, \nu, \mu}(t)\leq y]$ of the $Y_{\pi}^{\alpha, \nu, \mu}(t)$ is given by
        \begin{align*}
        Q_{\pi}^{\alpha, \nu, \mu}(t) &= \frac{1}{1+\mu t(1-e^{-(\lambda+\nu)^{\alpha}+\nu^\alpha})}\mathbb{I}_{y \geq 1}\\
        &+\frac{1}{1+\mu t(1-e^{-\nu^\alpha})}\sum_{m=1}^{\infty}(-\lambda)^m \mathbb{P}[V_k \leq y]\sum_{r=0}^{\infty}{\alpha r \choose m }\frac{(-(\lambda+\nu)^\alpha)^r}{r!} w_{r}\left(\frac{\mu t e^{-\nu^\alpha}}{1+\mu t(1-e^{-\nu^\alpha})}\right)
    \end{align*}
    \begin{proof}
    This is based on a similar approach as shown in \eqref{pr4.5}.
    \end{proof}
    \end{proposition}
    \begin{corollary}
        If $X_i's$ are discrete, then pmf $\mathbb{P}[Y_{\pi}^{\alpha, \nu, \mu}(t)=k]$ of the $Y_{\pi}^{\alpha, \nu, \mu}(t)$ is given 
        \begin{align}
            \mathbb{P}[Y_{\pi}^{\alpha, \nu, \mu}(t)=k]&= \frac{1}{1+\mu t(1-e^{-(\lambda+\nu)^{\alpha}+\nu^\alpha})}\mathbb{I}_{k= 1}\nonumber\\
        &+\frac{1}{1+\mu t(1-e^{-\nu^\alpha})}\sum_{m=1}^{\infty}(-\lambda)^m \mathbb{P}[V_m=k]\sum_{r=0}^{\infty}{\alpha r \choose m }\frac{(-(\lambda+\nu)^\alpha)^r}{r!} w_{r}\left(\frac{\mu t e^{-\nu^\alpha}}{1+\mu t(1-e^{-\nu^\alpha})}\right)\nonumber
        \end{align}
    \end{corollary}
    \begin{corollary}Let $\phi_{V_m}$ be the density of the $V_m$, then
    \begin{itemize}
        \item For absolutely continuous case, the density $h_{\pi}^{\alpha,\nu, \mu}(y, t)$ of $Y_{\pi}^{\alpha, \nu, \mu}(t)$ such as
        $$
        h_{\pi}^{\alpha, \nu, \mu}(y, t)=\frac{1}{1+\mu t(1-e^{-\nu^\alpha})}\sum_{m=1}^{\infty}(-\lambda)^m \phi_{V_m}(y)\sum_{r=0}^{\infty}{\alpha r \choose m }\frac{(-(\lambda+\nu)^\alpha)^r}{r!} w_{r}\left(\frac{\mu t e^{-\nu^\alpha}}{1+\mu t(1-e^{-\nu^\alpha})}\right),\;\; y \neq 0, \; t >0.
        $$
        \item When $Y_{\pi}^{\alpha, \nu, \mu}(t)=1$ then
        $$\mathbb{P}[Y_{\pi}^{\alpha, \nu, \mu}(t)=1]=\frac{1}{1+\mu t(1-e^{-(\lambda+\nu)^{\alpha}+\nu^\alpha})}.
 $$
        \end{itemize}
         \end{corollary}

\section{Applications in shock models}\label{sec:appl}

This section highlights the practical application of previously derived results in the context of shock models, which serve as essential tools for analyzing the behaviour of system lifetimes in random environments. Shock models are particularly relevant for understanding how systems endure and fail under harmful events. Two critical aspects form the foundation of any shock modelling approach. One is pattern of shock arrivals: this refers to the stochastic process governing when shocks occur. It could be modelled using simple frameworks such as Poisson processes, or renewal processes. Second is effect of shocks on the system: this captures the system's response to each shock, varying from negligible to catastrophic.

Based on how shocks impact the system, several classes of shock models have been extensively studied in the literature. These models can be broadly categorized into the following five types. $(i)$ Extreme shock models: in these models, the system fails immediately if a shock exceeds a certain critical threshold, regardless of prior shocks. $(ii)$ Cumulative shock models: here, the cumulative effect of successive shocks determines the system's failure. A system fails when the total accumulated damage surpasses its endurance limit. $(iii)$ Run shock models: in these model, the system fails when $k$ number of successive shocks exceeds the prefixed threshold value. $(iv)$ $\delta$ shock models: according to these models, the system fails when time lag between two consecutive shocks are too close or too far. Lastly, $(v)$ mixed shock models: these combine characteristics from two or more categories above, enabling greater flexibility to describe complex system behaviours. Some recent contributions can be found in Cha and Finkelstein~\cite{cf}, Chadjiconstantinidis and Eryilmaz~\cite{chad3},  Farhadian and Jafari~\cite{Farhadian}, Goyal \emph{et al.}~\cite{gxg}, Ozkut~\cite{o}, Soni \emph{et al.}~\cite{sp}, and references therein.
We refer Nakagawa~\cite{n}, and Cha and Finkelstein~\cite{cf1} as a good source of knowledge on this topic.

In the literature most shock models have been studied under the assumption of the Poisson process; however, this process has its own limitations. For instance, this process has independent increment property and stationary increment property. To overcome these limitations various generalizations of the Poisson process have been considered in the literature. For example, to capture correlation between the increments, researchers modelled shock processes through mixed Poisson processes (see, e.g., Cha and Finkelstein~\cite{cf},  Syuhada \emph{et al.}~\cite{sth}, Goyal \emph{et al.}~\cite{gfh}, to name a few). One reason of considering these processes is that they have positive dependency in its increments and another is that they provide mathematical tractable results. Recently, Goyal \emph{et al.}~\cite{gxg} have studied a class of shock models under the assumption of Markovian arrival process. 

Due to complexity in natural situations/or nature it is always interest of researchers to model shocks arrival pattern through a suitable counting process that fit well in the required natural situation. For example, as literature pointed out, there are many situations where non-fractional process models fits weaker than fractional processes. Therefore, it is natural interest of researchers of modelling shock arrivals using fractional counting processes. To best of our knowledge, Goyal \emph{et al.}~\cite{ghf} were first to study shock models based on fractional renewal process. In this study they study a class of shock models which contains many known shock models and provide an application to optimal maintenance by considering optimal age replacement policy. In this section, we are moving one step forward in this direction. Here we consider a general fractional counting process, defined in previous section, which contains several other known fractional processes as a special cases. In this section we use results of previous section to study two shock models; namely, extreme shock models and cumulative shock models.

Now we provide a brief description of the system under consideration. Assume that the system under consideration is working under harmful environment and can fail at any unexpected time. Shocks (external) are the only reason for system failure. Let $S$ denote the lifetime of the system. Let $Z_{i}$ denote the magnitude of the $i$th shock and $W_{i}$ represent damage given by the $i$th shock, $i \in \mathbb{N}$. Let the random variable $M(t)$ counts the number of shock in the time interval $(0, t]$, $t \geq 0$. Assume that the sequence $\{Y_{i}: i \in \mathbb{N} \}$ and $\{W_{i}: i \in \mathbb{N} \}$ contain independent and identically distributed random variables. Further assume that $\{Y_{i}: i \in \mathbb{N} \}$ and $\{W_{i}: i \in \mathbb{N} \}$ are independent from the shock process $\{M(t): t \geq 0\}$, respectively. With these assumptions, we first study the extreme shock model and then the cumulative shock model.
\subsection{Extreme shock model}
The extreme shock model has received considerable attention in reliability studies. In this model, the system fails if an individual shock exceeds a predefined threshold, such as a vehicle axle failing when a crack reaches a critical depth. That is, when the magnitude of $i$th shock exceed the predetermined threshold $\gamma$, or $\mathbb{P}(Y_{i} > \gamma)$ the system fails at $i$th shock, $i \in \mathbb{N}$. In the following theorem, we provide expression of the reliability function when the shock process follows the generalized counting process with the GCP.

%That is, if the magnitude of $i$th shock exceed the predetermined threshold $\gamma$ the system fails at $i$th shock, $i$ in $\mathbb{N}$.

\begin{theorem}
Suppose that the shock process $\{M(t): t \geq 0\}$ follows the generalized counting process with the GCP with the set parameters $\{f(\cdot), \lambda, \mu \}$. Then the reliability function of the lifetime of the system is given by:
$$\mathbb{P}(S > t) =  \frac{1}{1 + \mu t \left(1 - e^{-f(\lambda(1- q))} \right) }, $$
where $q$ represents the probability of survival of the system from a shock. In addition, the corresponding failure rate can be expressed as:
$$r_{S}(t) =  \frac{\mu \left(1 - e^{-f(\lambda(1- q))} \right) }{1 + \mu t \left(1 - e^{-f(\lambda(1- q))} \right) }. $$
\end{theorem}
\begin{proof}
Note that, according to the extreme shock model, probability of the system survival at time $t$ given that $k$ shocks has arrived by time $t$ can be written as
$$\mathbb{P}(S > t|M(t) = k) = q^{k}, \quad k = 0,1,2,\dots,$$
where $q$ is the probability of survival of each shock. Consequently,
$$ \mathbb{P}(S > t) =  \sum_{k=0}^{\infty} q^{k} P(M(t) = k).$$
From Equation~\eqref{deq1}, we get
\begin{eqnarray}\label{deq2}
\mathbb{P}(S > t) &=& \sum_{k=0}^{\infty}  \sum_{n=0}^{\infty} q^{k} \frac{1}{1+\mu t} \frac{(-1)^k}{k!}\frac{d^k}{du^k}e^{-nf(\lambda u)}|_{(u=1)} \left(\frac{\mu t}{1+\mu t}\right)^n  \nonumber \\
&=& \sum_{n=0}^{\infty}  \frac{1}{1+\mu t}   \left(\frac{\mu t}{1+\mu t}\right)^n  \left(\sum_{k=0}^{\infty}  q^{k} \frac{(-1)^k}{k!}\frac{d^k}{du^k}e^{-nf(\lambda u)}|_{(u=1)}\right).
\end{eqnarray}
Observe that, from the Taylor's series expansion, the expression
$$\sum_{k=0}^{\infty}  q^{k} \frac{(-1)^k}{k!}\frac{d^k}{du^k}e^{-nf(\lambda u)}|_{(u=1)} = e^{-nf(\lambda(1-q))}$$
holds. Hence, from Equation~\eqref{deq2}, we can write
$$\mathbb{P}(S > t) = \sum_{n=0}^{\infty}  \frac{1}{1+\mu t}   \left(\frac{\mu t}{1+\mu t}\right)^n e^{-nf(\lambda(1-q))}.$$
After simplifying the above expression, we get the required result for $\mathbb{P}(S > t)$. The expression of the failure rate follows directly from the expression of $\mathbb{P}(S > t).$
\end{proof}

In the following remark, we provide an alternative modeling to study the extreme shock model. The objective is to include this remark in order to show the applicability of the multiplicative process defined in the paper in shock modeling to study system's lifetime behavior. 

\begin{remark}
    Assume that the system under consideration has only two working states; namely 0 and 1. Here the 0 state represents that the system is failed while the state 1 represents that the system is under working condition. Let $\Phi(t)$ be the state of the system at any time $t \geq 0$. Then obviously, $\Phi(t) = 0 \;or\; 1$ for $t \geq 0$. Suppose that the system is subjected to shocks, and they can occur on the system at any random time. Any shock on the system either fail the system with probability $p$ or does not effect it with the compliment probability $q  = 1 - p$. Let $X_{i}$ be the state of the system at $i$th shock, $i \in \mathbb{N}$. Hence, $\mathbb{P}(X_{i} = 1) = q$ and $\mathbb{P}(X_{i} = 0) = p$, $i \in \mathbb{N}$. Therefore, the relation
    $$\Phi(t) = \prod_{i=1}^{M(t)} X_{i},$$
    holds for any $t \geq 0$. Now if $S$ represents the system lifetime then the reliability function of $L$ is given by
    $$\mathbb{P}(S > t) = \mathbb{P}(\Phi(t) = 1) = \mathbb{E}(\Phi(t)).$$ Suppose Suppose that the shock process $\{M(t): t \geq 0\}$ follows the generalized counting process with the GCP with the set parameters $\{f(\cdot), \lambda, \mu \}$. Then From~\eqref{deqmcp}, we get
    $$\mathbb{P}(S > t) =  \frac{1}{1 + \mu t \left(1 - e^{-f(\lambda(1- q))} \right) }, $$
    because $\mathbb{E}(X_{1}) = q$.  
\end{remark}

\begin{figure}[h!]
 \centering
\begin{subfigure}{.5\textwidth}
\includegraphics[width=\linewidth]{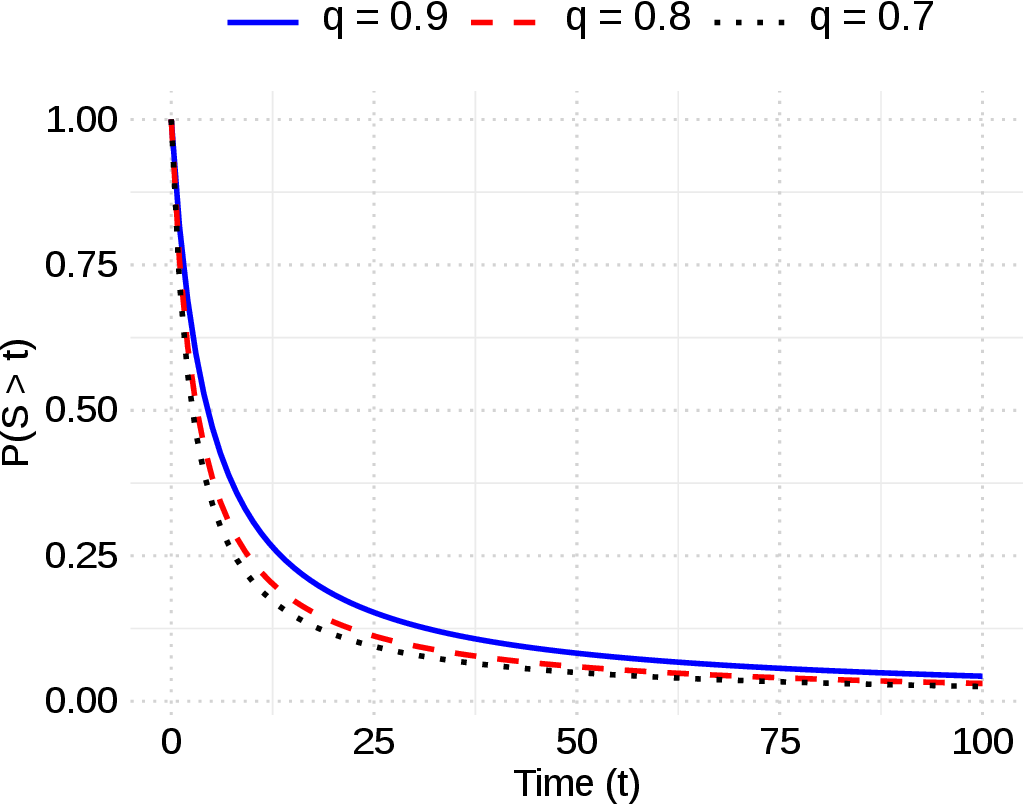}
\caption{ for $q$}
\label{fig11}
\end{subfigure}%
\begin{subfigure}{.5\textwidth}
\includegraphics[width=\linewidth]{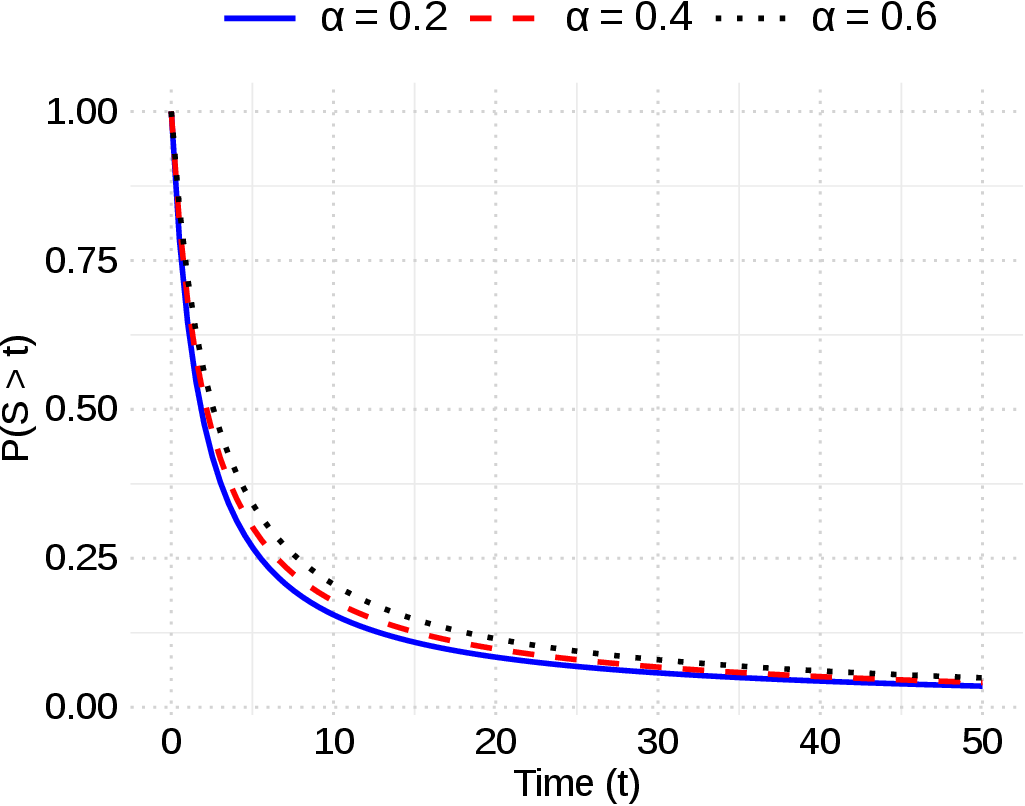}
  \caption{ for $\alpha$}
  \label{fig12}
\end{subfigure}
\begin{subfigure}{.5\textwidth}
\vspace{0.5cm} % Adjust the vertical space here
\includegraphics[width=\linewidth]{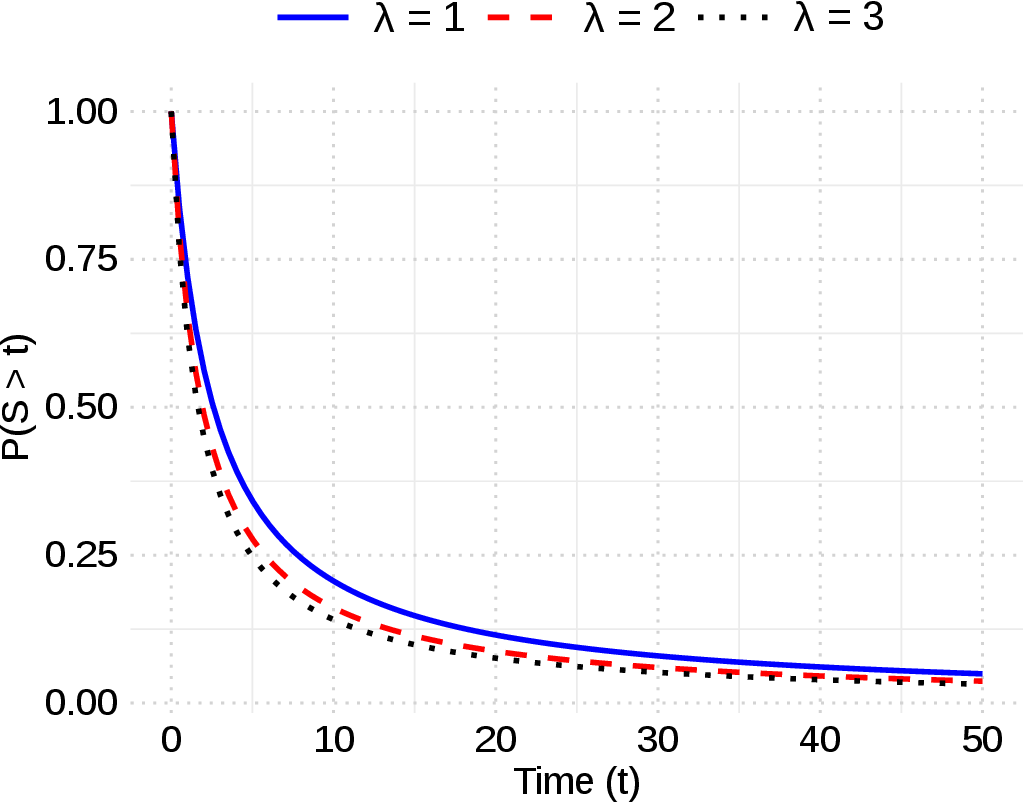}
\caption{ for $\lambda$}
\label{fig13}
\end{subfigure}%
\begin{subfigure}{.5\textwidth}
\includegraphics[width=\linewidth]{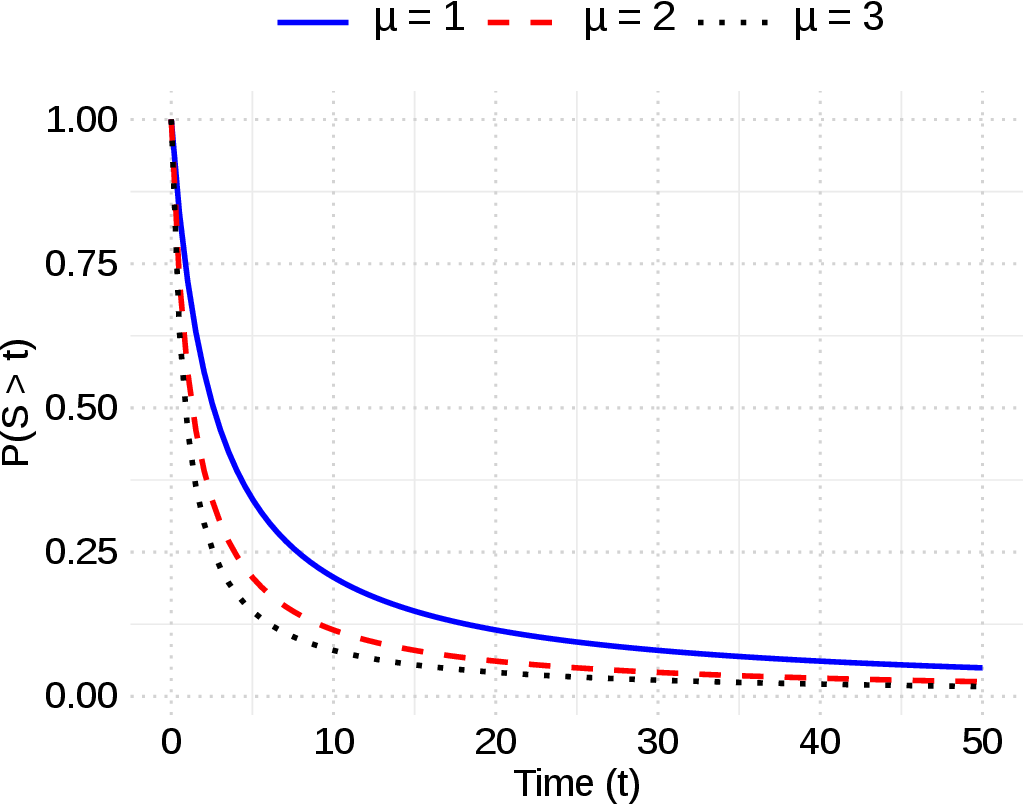}
  \caption{ for $\mu$}
  \label{fig14}
\end{subfigure}
\caption{Sensitivity analysis of the system reliability for the case when $f(s) = s^{\alpha}$, $0 < \alpha < 1$.}
\label{fig1}
\end{figure}

\begin{figure}[h!]
 \centering
\begin{subfigure}{.5\textwidth}
\includegraphics[width=\linewidth]{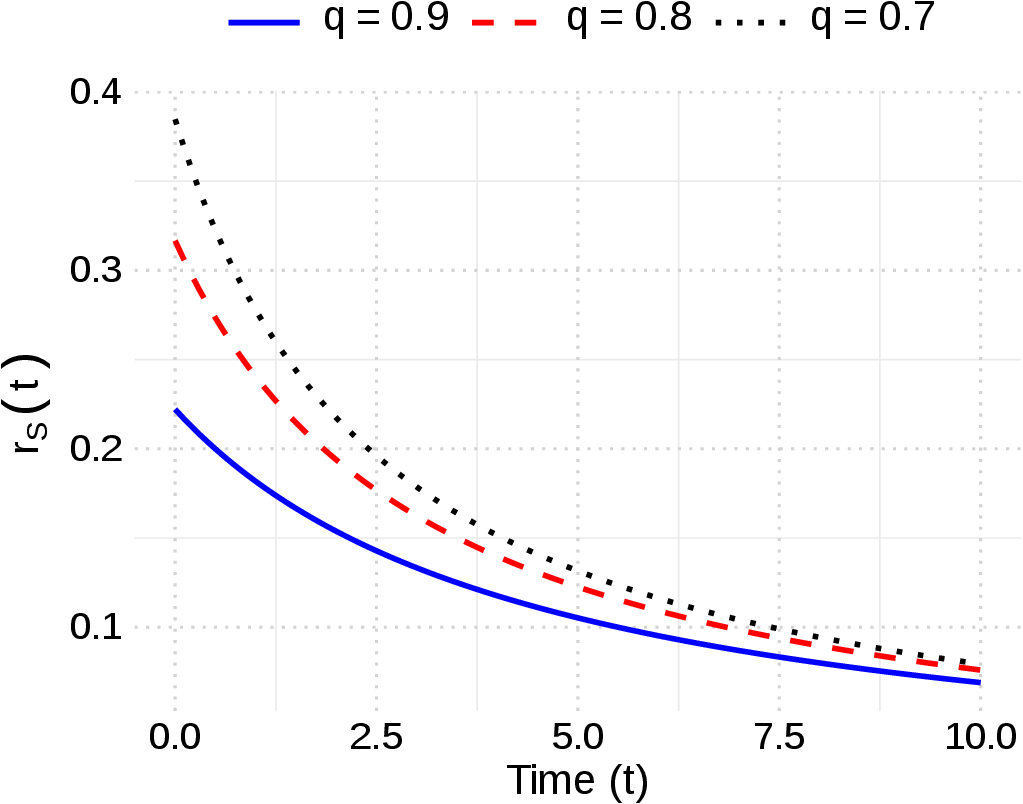}
\caption{ for $q$}
\label{fig21}
\end{subfigure}%
\begin{subfigure}{.5\textwidth}
\includegraphics[width=\linewidth]{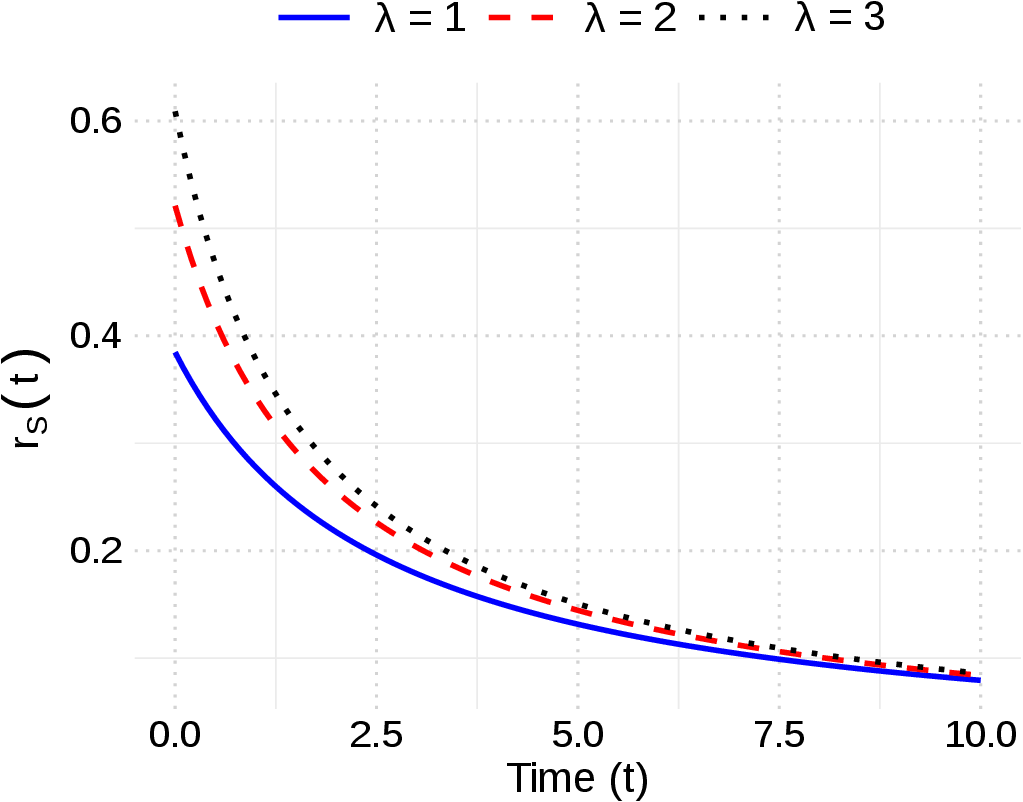}
  \caption{ for $\alpha$}
  \label{fig22}
\end{subfigure}
\begin{subfigure}{.5\textwidth}
\includegraphics[width=\linewidth]{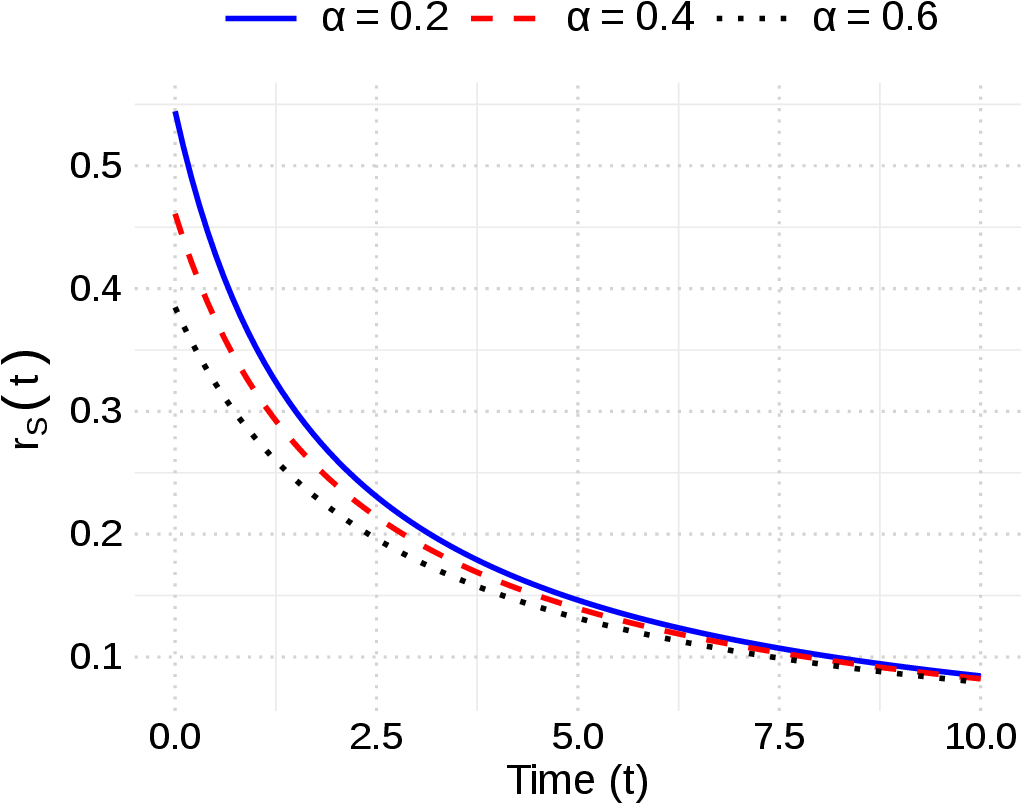}
\caption{ for $\lambda$}
\label{fig23}
\end{subfigure}%
\begin{subfigure}{.5\textwidth}
\includegraphics[width=\linewidth]{failure_plot_wrtal.eps}
  \caption{ for $\mu$}
  \label{fig24}
\end{subfigure}
\caption{Sensitivity analysis of the system's failure rate for the case when $f(s) = s^{\alpha}$, $0 < \alpha < 1$.}
\label{fig2}
\end{figure}

Notice that the above result generally hold for any choice of function $f(\cdot)$. Hence survival function of system lifetime can be derived for special cases given in~\eqref{Levy_exponent}. For demonstration of results we consider $f(s) = s^{\alpha}$ (stable subordinator), where $0 < \alpha < 1$.  However, results can be demonstrate for other cases as well. We illustrate our result in the Figure~\ref{fig1} and Figure~\ref{fig2}. From Figure~\ref{fig11}, we can see that decreasing the value of the parameter $q$ decrease the system reliability; this is because $q$ denotes the probability of system survival from a shock, as a result larger value of $q$ enhance system survivability. Figure~\ref{fig12} shows that increment in values of $\alpha$ increase the system reliability. Further, Figures~\ref{fig13} and~\ref{fig14} shows that increment in the value of either $\lambda$ or $\mu$ decrease the the system reliability. Figure~\ref{fig2} shows sensitivity analysis of the system's failure rate; similar conclusion as reliability function can be drawn from this figure. For plotting Figures~\ref{fig1} and~\ref{fig2}, we consider the same baseline parameters: $q = 0.7$, $\alpha = 0.6$, $\lambda = 1$, and $\mu = 1$. We perform the analysis independently for each parameter while holding the others fixed.

\subsection{Cumulative shock model}
In cumulative shock models, system failure occurs when the accumulated damage from successive shocks surpasses a predefined threshold $T$. These models are particularly relevant in various real-world applications. For example, the strength of a fibrous carbon composite depends significantly on the strength of its individual fibers, which may break under tensile stress. The composite material ultimately fails as a result of cumulative damage (see Nakagawa~\cite{n}, p. 2). 
Let $Z(t)$ represent the total cumulative damage at time $t$. According to the cumulative shock model, this total damage is described as a random sum:
$$
Z(t) = \sum_{i=1}^{M(t)} W_i,
$$
where it is conventionally assumed that $ \sum_{i=1}^{0} Y_i = 0 $.

In this section, we consider a cumulative shock model where damage size is discrete with set of non-negative integers as a support. Such kind of model may useful in many real-world applications. Some examples are as follows:
\begin{enumerate}
    \item[$(i)$]  In a multi-component system, a shock can fail components of the system. Thus, in this case, number of failed component can be consider damage size for a particular shock. If total number of component failure is more than a preset threshold number then the multi-component system can fail.
    \item[$(ii)$] Electric vehicles battery packs are composed of numerous individual cells arranged in series and parallel configurations. Each cell contributes a specific amount of voltage and capacity to the overall battery pack. Over the time, individual cells may fail completely due to some external factors such as: thermal stress, mechanical stress, etc. When a cell fails, it reduces the overall capacity and efficiency of the battery pack. When the battery pack cross its tolerable limit threshold, the pack may be considered failed and require repair or replacement.

\end{enumerate}

In the following theorem, we provide expression of the reliability function of the system under cumulative shock model.

\begin{theorem}\label{dthmcm}
Let us assume that $\{M(t): t \geq 0\}$ is a geometric process with parameters $\mu$ and $\lambda$. Further, assume that damage size $W_{i}$ follows subordinated Poisson random variables with the parameter $f(\cdot)$. Suppose $T$ is the threshold of the maximum damage that system can tolerate. Then the reliability function of $S$ is given by
$$ \mathbb{P}(S > t) = \sum_{k=0}^{T-1}  \frac{(-1)^k}{k!}\frac{d^k}{du^k}\left(\frac{1}{1+\mu t(1-e^{-f(\lambda u)})}\right)_{u=1}.$$
\end{theorem}
\begin{proof}
    First note that, from the definition of the cumulative shock model, reliability function of the system lifetime and the cumulative damage has the following relation:
    \begin{eqnarray}\label{dthm2}
        \mathbb{P}(S >t) = \mathbb{P}\left(Z(t) < T\right).
    \end{eqnarray}
    Now from Propositions~\ref{pr1} and~\ref{pr2}, we get
    $$ \mathbb{P}(Z(t) = k) = \frac{(-1)^k}{k!}\frac{d^k}{du^k}\left(\frac{1}{1+\mu t(1-e^{-f(\lambda u)})}\right)_{u=1}.$$
    Thus, result follows directly from Equation~\eqref{dthm2}.
\end{proof}

\ifx

To illustrate the result given in the above theorem, we consider a special case $f(s) = s^{\alpha}$, $0 < \alpha <1$. For this special case, by Proposition~\ref{pmf_sfpp_gcp} and Theorem~\ref{dthmcm}, system's reliability can be expressed as:
 $$
  \mathbb{P}(S > t) = \sum_{k=0}^{T-1} \sum_{r=0}^\infty {\alpha r \choose k} \frac{(-\lambda^\alpha)^r (-1)^k }{r!}w_r(\mu t),\;\; \lambda>0, \; \mu >0.
  $$
  Here $w_{r}(\cdot)$ is the geometric polynomial of degree $r$.

  \fi

% \newpage
%\bibliographystyle{abbrv}
%\bibliography{researchbib}
%\clearpage

\def\cprime{$'$}

\end{document}